\documentclass[a4paper,twocolumn]{autart}

\usepackage{calrsfs,times}
\usepackage{graphicx,enumerate}
\usepackage{psfrag}
\usepackage{pstool}
\usepackage{amssymb}
\usepackage{amsmath}
\usepackage{calrsfs}
\usepackage{theorem}
\usepackage{color}
\usepackage{epsfig,epstopdf,wrapfig}
\usepackage{paralist}
\usepackage{float,cuted}
\usepackage{multirow}
\usepackage{subcaption}

\newtheorem{definition}{Definition}
\newtheorem{theorem}{Theorem}
\newtheorem{lemma}{Lemma} 
\newtheorem{corollary}{Corollary} 
\newtheorem{assumption}{Assumption}
\newtheorem{proposition}{Proposition}

{\theorembodyfont{\upshape}}
{\theorembodyfont{\upshape}\newtheorem{remark}{Remark}}
{\theorembodyfont{\upshape}}
{\theorembodyfont{\upshape}\newtheorem{algo}{Algorithm}}

\usepackage{amsmath}
\usepackage{accents}
\newlength{\dhatheight}

\usepackage{color}
\usepackage[normalem]{ulem}
\usepackage{soul} 
\soulregister\cite7
\soulregister\ref7
\soulregister\pageref7

\begin{document}
\begin{frontmatter}

\title{
  Design of Coherent Passive Quantum Equalizers Using Robust Control Theory}
\thanks[footnoteinfo]{This work was 
    supported by the Australian Research Council under the Discovery Projects funding scheme (project DP200102945).}

\author[First]{V.~Ugrinovskii}\ead{v.ougrinovski@adfa.edu.au} \and
\author[Second]{M.~R.~James}\ead{Matthew.James@anu.edu.au}
\address[First]{School of Engineering and Technology, University of New
  South Wales Canberra, Canberra, ACT, 2600, Australia}
\address[Second]{School of Engineering, College of Engineering, Computing
  and Cybernetics, The Australian National University, Canberra, ACT 2601,
  Australia}

\begin{abstract}
The paper develops a methodology for the design of coherent
equalizing filters for quantum communication channels. Given a linear quantum
system model of a quantum communication channel, the aim is to obtain another
quantum system which, when coupled with the original system, 
mitigates degrading effects of the environment. The main result of the
paper is a systematic equalizer synthesis algorithm which
relies on methods of state-space robust control design via semidefinite
programming.   
\end{abstract}

\end{frontmatter}

\section{Introduction}

\begin{figure}[t]
\begin{center}
\psfragfig[width=0.8\columnwidth]{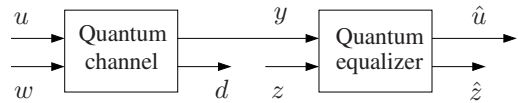}{ 
\psfrag{Quantum}{\hspace{-1.5ex}\footnotesize Quantum}
\psfrag{channel}{\hspace{-1.5ex}\footnotesize channel}
\psfrag{equalizer}{\hspace{-1.5ex}\footnotesize equalizer}
  \psfrag{+}{$+$}
  \psfrag{-}{$-$}
  \psfrag{b}{$u$}
  \psfrag{w}{$w$}
  \psfrag{what}{$d$}
  \psfrag{z}{$z$}
  \psfrag{bhat}{$\hat{u}$}
  \psfrag{zhat}{$\hat{z}$}
  \psfrag{y}{$y$}}
  \caption{A general quantum communication system.
  } 
  \label{fig:general}
\end{center}
\end{figure}

This paper develops a robust control approach to derivation of
coherent equalizers for quantum communication channels. The class of
systems we consider 
includes a wide range of linear quantum systems comprised of
quantum optical components such as beamsplitters, optical cavities, phase
shifters, that may be used in quantum communication systems. As messages
encoded in a beam of light are transmitted through such 
systems, they degrade due to losses to the environment and distortions
caused by the channel itself. Coherent equalization seeks to
mitigate environmental and channel
distortions blue in the system by coupling it with another
quantum physical system acting as a 
filter. A general diagram of a quantum equalization system is shown in
Fig.~\ref{fig:general}. As we explain below, in this figure the symbols $u,w$
(respectively, $y,z$) represent the input 
field and the environment of the channel (respectively, the equalizing
filter), $d$, $\hat z$ represent losses to the system and filter
environments, and $\hat u$ represents the output field of the filter.  
The aim is to obtain a filtering device which minimizes the
  mismatch between $u$ and $\hat u$ in the mean-square sense. 

From the system theoretic viewpoint, the above task 
is analogous to the task of equalization of classical communication
channels~\cite{HSK-1999}, with a 
notable difference that the filtering is done at the physical level using a
quantum device rather than a classical (i.e., nonquantum) signal processor. 
Nevertheless, the analogy prompts the question whether the powerful
Wiener's mean-square optimization
paradigm~\cite{Kailath-1981,KSH-2000,Wiener-1949} 
widely used in classical communications can 
be extended into the realm of nonclassical communication quantum systems
design.   
This problem has been introduced recently in~\cite{UJ2b,UJ2a,UJ2} in the
context of optimization of the mean-square mismatch between
the system input and the filter output. We also mention earlier results in a
similar vein which were concerned with developing coherent versions of the
Luenberger observer and Kalman filter~\cite{MJP-2016,VP-2013,VA-2014}. 
However these developments are significantly different from the coherent
equalization problem introduced in~\cite{UJ2b,UJ2a,UJ2} and also considered
in this paper in that they are aimed at estimation of the internal 
modes of the quantum channel system. In contrast, the coherent
equalization is concerned with optimally matching the input 
fields of the quantum system in the mean-square sense. The methods
in~\cite{MJP-2016,VP-2013,VA-2014} are not applicable to this
problem.  

To introduce the coherent equalization problem, consider a
linear quantum system consisting of a 
collection of quantum harmonic oscillators driven by quantum noise
fields representing the input to the system~\cite{GJN-2010}. We consider
\emph{completely passive} quantum systems, i.e., systems whose Hamiltonian and
coupling with the input fields involve only annihilation operators of the
oscillator modes~\cite{MP-2011,NY-2017}. This
terminology reflects the fact
that physical implementation of such systems does not require external
sources of quanta and can be done using only passive optical components 
such as beam splitters, phase shifters and mirrors~\cite{NY-2017}. The dynamic 
behaviour of such systems in the 
Heisenberg picture of quantum mechanics can be described 
by the quantum stochastic differential equation in the Langevin
form~\cite{GJN-2010,GZ-2015,MP-2011,NY-2017,WM-2008,WM-2009} 
\begin{eqnarray}
  \label{dyn}
  \dot{\mathbf{a}}(t)&=&A \mathbf{a}(t)+
                         B_1u(t)+B_2w(t), \nonumber \\ 
  y(t)&=&C_1 \mathbf{a}(t)+ D_{11}u(t)+ D_{12}w(t), \nonumber \\ 
  d(t)&=&C_2 \mathbf{a}(t)+ D_{21}u(t)+ D_{22}w(t).
\end{eqnarray}
Here $A$, $B_j$, $C_j$ and $D_{kj}$ are complex matrices, $j,k=1,2$, and
$\mathbf{a}$ is the vector of annihilation operators of the 
oscillator modes. Also, $u,w$ are vectors of quantum noise processes
corresponding 
to annihilation operators of the input field in an infinite-dimensional
Hilbert space called the Fock space~\cite{NY-2017}; $u$ represents the input fields 
engineered to carry the transmitted information, and $w$ is associated with the
physical environment of the system; see
Fig.~\ref{fig:general}. Furthermore, $y$, $d$ are the vectors 
of quantum noise processes corresponding to annihilation operators of the output
field of the system; $y$ represents the part of the output field from which
the transmitted information will be recovered while 
$d$ represents the loss to the environment.  
The transfer function $G(s)$ of the system~(\ref{dyn}) 
is square and paraunitary; 
it relates the bilateral Laplace transforms of 
$y(t)$, $d(t)$ and $u(t)$, $w(t)$~\cite{GJN-2010,ZJ-2013}. Therefore it can
be partitioned as 
\[
G(s)=
   \left[
     \begin{array}{cc}
 G_{11}(s) & G_{12}(s)\\ G_{21}(s) & G_{22}(s)
     \end{array}
   \right],
\] 
where $G_{jk}(s)=C_j(sI-A)^{-1}B_k+D_{jk}$, $j,k=1,2$.

A coherent filter is another linear quantum system coupled with the
designated output fields of the system~(\ref{dyn}) and its own
environment (see Fig.~\ref{fig:general}), with the aim of extracting information from the former. Such
system can be designed in a number of ways, and this paper follows the 
seminal Wiener's approach~\cite{Wiener-1949}. Namely, we are
concerned with 
derivation of a completely passive quantum system governed by
the output fields $y$ of the system~(\ref{dyn}) and its own environment $z$
which optimizes the power spectrum density (PSD) of the difference between the
field $u$ and the equalizer's output field $\hat u$ in the following sense:
\begin{equation}
    \sup_\omega \boldsymbol{\sigma} \left[\int_{-\infty}^\infty \langle
    e(t)e(0)^\dagger\rangle e^{-i\omega
    t}dt\right] \to \inf;
\label{eq:53}
\end{equation}
here $e(t)=\hat u(t)-u(t)$, $\boldsymbol{\sigma}[\cdot]$ is the largest
eigenvalue of a Hermitian matrix, and $\langle\cdot\rangle$ is
the operator of quantum expectation of the system in a Gaussian state
$\boldsymbol{\rho}$ evolving in a stationary regime; see~\cite{UJ2a,UJ2}
and Definition~\ref{def.error} below. 
The infimum is taken over the class of transfer functions of completely
passive quantum systems of interest. Formally, this class consists of     
rational transfer functions $H(s)$ of \emph{physically realizable} causal stable
linear mappings  
$\mathrm{col}(y,z)\to \mathrm{col}(\hat u,\hat z)$ (see Definition~\ref{D.1} below), where $y$ is the output
vector of the system~(\ref{dyn}) described above, and  
the vector quantum processes $z$ and $\hat u$, $\hat z$ correspond to the
annihilation operators of the filter environment and the filter output
fields, respectively, with $\hat u$ having the same dimension as $u$. Also, 
$\mathrm{col}(\cdot,\cdot)$ denotes concatenation of column vectors.
That is, we pose the optimization
objective in the frequency domain (since the PSD is a frequency function), in
contrast to,~e.g., \cite{VA-2014,VP-2013} where optimization objectives
were formulated in the time domain.   
 
Our main result is a tractable algorithm that allows to compute a
near optimal solution to the optimization problem~(\ref{eq:53}). Combined
with the existing implementation
results for completely passive quantum
systems~\cite{Nurdin-2010,NJD-2009,NY-2017}, this provides a  
systematic method for the design of coherent equalizers
for linear quantum systems. 

Realizability of the filter as a physical quantum system requires that only
causal linear mappings which preserve the  
canonical commutation relations of the underlying annihilation operators
and their adjoint operators  are
admissible~\cite{JNP-2008}. This 
requires the transfer function $H(s)$ in question to be
stable, causal and paraunitary (lossless bounded
real)~\cite{MP-2011}. These properties are the frequency domain 
design constraints on $H(s)$ as a decision variable in the above mentioned
optimization problem~\cite{UJ2a,UJ2b,UJ2}. They require $H(s)$ to
belong to the Hardy space $H_\infty$ and be paraunitary,  
\begin{equation}
  \label{eq:60}
H(s)^H H(s)=H(s)H(s)^H=I.
\end{equation}
Here and thereafter, $H(s)^H\triangleq H(-s^*)^\dagger$,
where the symbols~$^*$ and~$^\dagger$ denote  the complex conjugate number
and the complex conjugate transpose matrix, respectively. 
Other coherent filter and control design problems also account for
physical realizability
constraints~\cite{MJP-2016,VP-2013,VA-2014}. However, unlike this paper 
and~\cite{UJ2a,UJ2b,UJ2}, physical realizability is expressed in these
references as a set of nonconvex algebraic constraints on the matrices of
the filter/controller state space model.   

Since the physical realizability
constraint~(\ref{eq:60}) is not convex in general, the method for
minimizing the power spectrum density
of the equalization error proposed in this paper employs a convex
relaxation. It reduces the underlying optimization problem to  
an auxiliary problem involving minimization over a convex subset of the set of
strictly contractive (in the $H_\infty$ sense) causal mappings $y\to \hat
u$. The proposed relaxation has two significant features. Firstly, every
feasible transfer function of the auxiliary problem, if it exists, is shown
to lead to a 
physically realizable completely passive quantum system. The relaxation
considered in our previous work~\cite{UJ2} did not generally guarantee
this --- in order to give rise to a physically realizable $H(s)$, solutions to
that relaxed problem had to meet some additional conditions which had to
be checked after a solution had been found. Restricting 
the feasible set of the auxiliary problem to include only strictly contractive
mappings ensures that these conditions are satisfied whenever the auxiliary
problem has a solution, even though this is accomplished at the expense of
slightly reducing the class of candidate filters, compared with~\cite{UJ2}. 

The second feature of the proposed relaxation is that the filter transfer
functions constructed from the auxiliary problem provide a
quantitative upper bound on the theoretically smallest possible mean-square
equalization accuracy that can be achieved within the class of completely passive
coherent filters, i.e., 
they are bona fide \emph{guaranteed cost} filters. Moreover, we show that
under certain conditions this upper bound is exact in the sense that
filters can be constructed based on the proposed relaxation which
approximate the theoretically optimal equalization accuracy arbitrarily
closely. 
The systems which satisfy these conditions include quantum communication
systems~(\ref{dyn}) in which the intensity of the environment noise is
sufficiently high, and respectively, the signal-to-noise ratio (SNR) is
sufficiently low. 
This is illustrated in Section~\ref{examples}, which demonstrates
application of the proposed method to equalization of
quantum optical systems consisting of optical cavities and
beam splitters.

The key technical assumption which underpins our approach to the filter design
requires a certain transfer function matrix
associated with the system~(\ref{dyn}) to admit a spectral
factorization; see Assumption~\ref{A.rho}. The
counterpart of this technical assumption in~\cite{UJ2} 
required a similar transfer function matrix to admit a $J$-spectral
factorization. Thus, the classes of quantum systems captured by the two
approaches are not identical. 
For instance, the size of the vector $y$ in this paper does not have to be
the same as the size of $u$, whereas in~\cite{UJ2} these vectors were to
have the same dimension. Moreover, the factorization assumption adopted
here 
enables us to 
recast the auxiliary convex relaxation problem within the 
semidefinite programming framework, as a convex
optimization problem involving a finite number of decision variables  
and linear matrix inequality constraints~\cite{LMI}. Semidefinite programs
of this kind are akin to those arising in the $H_\infty$ control
design~\cite{DP-2000}, and powerful numerical tools exist to solve them
efficiently~\cite{LMI,Nesterov-Nemirovskii}. Combined with the existing
methods for computing spectral factors~\cite{Anderson-1967,CALM-1997}, this
leads to a systematic, numerically tractable algorithm for the synthesis of
coherent completely passive equalizing filters.     

The paper is organized as follows. The constrained optimization formulation
of the coherent equalization problem is  
presented in Section~\ref{sec:equal-probl}. The analysis of the auxiliary
optimization problem and a sufficient condition for suboptimal solutions to
this auxiliary problem to have the desired contractiveness properties are presented
in Section~\ref{relaxation}. Section~\ref{framework} derives a semidefinite
program for computing such a solution. The algorithm for the synthesis of
physically realizable near optimal equalizer transfer functions is
given in 
Section~\ref{X}. Section~\ref{examples} 
presents examples which illustrate the 
proposed synthesis method. The first example presents a complete analysis of a
quantum optical system consisting of an optical cavity and several beam
splitters. We also show how a suboptimal filter for the 
system in this example can be implemented using standard quantum optical
components. In the second example, an equalizer is obtained for an
interconnection of two such systems. Conclusions and suggestions for
future work are given in Section~\ref{Conclusions}.

\section{Coherent equalization problem for quantum linear 
  systems}\label{sec:equal-probl}

\subsection{An open quantum system model of a quantum communication channel}

Consider a completely passive quantum system described by stochastic
differential equation~(\ref{dyn}). The quantum noise processes $u$ and $w$
of the input field have dimensions $n$, $n_w$, respectively. These
operators are assumed to satisfy the canonical commutation relations 
$[u_j(t),u_k^*(t')]=\delta_{jk}\delta(t-t')$,
$[w_j(t),w_k^*(t')]=\delta_{jk}\delta(t-t')$;
here $\delta(t)$ is the Dirac delta function, and $\delta_{jk}$ is the
Kronecker symbol: $\delta_{jk}=1$ when $j=k$, otherwise $\delta_{jk}=0$.
Also, $[u_j(t),w_k^*(t')]=0$. That is, the signal and environment operators
commute.  
We also assume that the system is in a Gaussian thermal state, and $u$, $w$
are zero mean Gaussian fields\footnote{For convenience of presentation, we
  will not make distinction between quantum fields and the respective vectors
  of annihilation operators and quantum noise processes. We will refer to
  them generically as inputs or outputs where this does not lead to a
  confusion.}; 
i.e., $\langle 
u(t)\rangle=0$, $\langle w(t)\rangle=0$ where $\langle \cdot \rangle$ is
the quantum 
  expectation of the system in this
  state~\cite{Parthasarathy-2012}. Furthermore, it is assumed  
that the processes $u$ and $w$ are 
not correlated, $\langle u(t)w^\dagger
(t')\rangle=0$. Letting $u^\#$, $w^\#$ denote the column vectors comprised
of the adjoint operators of $u$, $w$, and letting $\breve
u=\mathrm{col}(u,u^\#)$, $\breve w=\mathrm{col}(w,w^\#)$ be the vectors of 
operators obtained by concatenating $u$, $u^\#$ and $w$, $w^\#$,
respectively,  introduce the autocorrelation matrices of the processes $\breve u$, $\breve w$: 
\begin{eqnarray}
  \label{eq:44}
&&  R_{\breve{u}}(t)\triangleq \langle
  \breve{u}(t)\breve{u}^\dagger(0)\rangle=\left[
  \begin{array}{cc}
   I+\Sigma_{u}^T & 0 \\
   0 & \Sigma_{u} 
  \end{array}
  \right]\delta(t), \nonumber \\
&&  R_{\breve{w}}(t)\triangleq \langle
  \breve{w}(t)\breve{w}^\dagger(0)\rangle=\left[
  \begin{array}{cc}
   I+\Sigma_{w}^T & 0 \\
   0 & \Sigma_{w} 
  \end{array}
  \right]\delta(t).
\end{eqnarray}
The Hermitian positive definite $n\times n$ and $n_w\times n_w$
matrices $\Sigma_u$, $\Sigma_w$ symbolize the intensity of the signal and
noise in the system~(\ref{dyn}), respectively.

The matrices $A$, 
$B\triangleq [B_1~B_2]$, 
\[
C\triangleq
\left[
  \begin{array}{c}
    C_1 \\ C_2
  \end{array}
\right], \quad 
D\triangleq
\left[
  \begin{array}{cc}
    D_{11} & D_{12} \\ D_{21} & D_{22}
  \end{array}
\right]
\]
of the system~(\ref{dyn}) are complex matrices of dimensions $m\times m$,
$m\times (n+n_w)$, $(n_y+n_d)\times m$ and $(n_y+n_d)\times (n+n_w)$
respectively. The row dimensions of $C_1$, $C_2$, denoted $n_y$,
$n_d$, and the corresponding row dimensions of the blocks of the matrix $D$
do not have to be equal to $n$, $n_w$, however it holds that
$n_y+n_d=n+n_w$, i.e., the total number of input fields is equal
to the total number of output fields; cf.~\cite{UJ2}. The
coefficients of the system~(\ref{dyn}) are assumed to satisfy equations
reflecting the physical realizability of the system,  
\begin{eqnarray}
  \label{eq:42}
  && A+A^\dagger+BB^\dagger=0, \nonumber \\
  && B=-C^\dagger D, \quad D^\dagger D=DD^\dagger =I_{n+n_w},
\end{eqnarray}
\cite{SP-2012,JNP-2008,MP-2011}; $I_k$ is the $k\times k$ identity
matrix. Without loss of generality, we assume that $(A,B,C,D)$ is the minimal
realization of the transfer function $G(s)\triangleq D+C(sI_m-A)^{-1}B$, and so the matrix $A$ is
Hurwitz~\cite{GZ-2015}. 

The physical realizability of the 
system~(\ref{dyn}) dictates that the transfer function $G(s)$ is
paraunitary, and $G(i\omega)$ is a unitary
matrix~\cite{SP-2012,GJN-2010,BGR-2013,MP-2011}:  
\begin{eqnarray}
  \label{eq:32}
  &&G(s)^H G(s)=G(s)G(s)^H=I_{n+n_w}, \\
  &&G(i\omega)^\dagger G(i\omega)=G(i\omega)G(i\omega)^\dagger=I_{n+n_w}.
  \label{eq:32.w}
\end{eqnarray}
Then $G(i\omega)$ is bounded at infinity and
analytic on the entire closed imaginary
axis~\cite[Lemma~2]{Youla-1961}. 

\subsection{Coherent equalization problem} 

As in~\cite{UJ2a,UJ2}, we consider completely passive filters for the
system~(\ref{dyn}). The filter environment is represented by a vector of
$n_z$ quantum noise operators $z$ and their corresponding adjoint operators,
with the canonical commutation relations 
$[z_j(t),z_k^*(t')]=\delta_{jk}\delta(t-t')$. For simplicity, we assume
that the environment is a zero mean noise in a Gaussian vacuum state. That
is, $\langle \breve z(t)\rangle =0$ where $\breve
z=\mathrm{col}(z,z^\#)$, and the correlation function of the noise
process $\breve z(t)$ is  
$  \langle
     \breve z(t) 
     \breve{z}^\dagger(t')\rangle=
  \left[
    \begin{array}{cc}
I_{n_z} & 0 \\
0 & 0
\end{array}\right]\delta(t-t')
$. 
It is assumed that $\breve z$ commutes with $\breve{u}$ and $\breve{w}$,
and $\langle \breve z(t) 
     \breve{u}^\dagger(t')\rangle=0$, $\langle \breve z(t) 
     \breve{w}^\dagger(t')\rangle=0$. 

The input into the equalizer combines the output fields of the
system~(\ref{dyn}) and the filter environment; its $n_f$-dimensional
vector of annihilation
operators is $\mathrm{col}(y,z)$, $n_f=n_y+n_z$. The vector of the filter
output processes has the same dimension $n_f$, but is
partitioned into vectors $\hat u$ and $\hat z$ so that  $\hat u$ has the
same dimension $n$ as $u$. We 
designate the first component of this partition, namely $\hat u$,
to serve as the estimate of $u$.   

\begin{definition}\label{D.1}
  An element $H(s)$ of the Hardy space $H_\infty$ is said to represent an
  admissible 
  completely passive physically realizable equalizer if $H(s)$ is a stable
  rational $n_f\times n_f$ transfer function, $n_f=n_y+n_z\ge n$,
  which is analytic in the right half-plane  $\mathrm{Re}s>-\tau$ ($\exists
  \tau>0$) and is paraunitary in the sense of~(\ref{eq:60}). The
  set of admissible equalizers will be denoted $\mathcal{H}_p$.   
\end{definition}

\begin{remark}\label{Rem1}
According to Definition~\ref{D.1}, admissible transfer functions are analytic and bounded in the closed
right half-plane of the complex plane, therefore each such transfer
function defines a causal system~\cite{DP-2000}. Formally we require
analiticity in an open half-plane $\mathrm{Re}s>-\tau$ ($\exists
\tau>0$) to ensure compatibility with the development in~\cite{UJ2}. 
The latter reference inherited this form of analiticity from~\cite{Youla-1961}
whose results were used in the derivation of physically realizable
equalizers.
\end{remark}

\begin{definition}[\cite{UJ2a,UJ2}]\label{def.error}
The equalization error of the quantum channel-filter system under
consideration is $e(t)=\hat u(t)-u(t)$, i.e., the difference between the
stationary $n$-dimensional operator valued processes $u$, $\hat u$ designated
as the channel input and the filter output. The power spectrum density
$P_e(s)$ of the equalization error is the bilateral Laplace transform of
the autocorrelation function $R_e(t)=\langle e(t)e(0)^\dagger\rangle$. 
\end{definition}

By definition, the power spectrum density
$P_e(s)$ is an $n\times n$ para-Hermitian transfer function
matrix. Its restriction to the imaginary axis yields a Hermitian nonnegative
definite matrix $P_e(i\omega)$. It was shown in~\cite{UJ2a,UJ2} that
\begin{equation}
  \label{eq:59}
P_e(s)=
      \left[
      \begin{array}{cc}
        H_{11}(s) & I
      \end{array}
      \right] \Phi(s)       \left[
      \begin{array}{c}
        H_{11}(s)^H \\ I
      \end{array}
      \right]. 
\end{equation}
Here $H_{11}(s)$ is the top-left block in the partition of
the filter transfer function $H(s)$ into blocks compatible with the
partitions of the filter input $\mathrm{col}(y,z)$ and its output
$\mathrm{col}(\hat u, \hat z)$:    
\begin{equation}
  \label{eq:98a}
  H(s)=
  \left[
    \begin{array}{cc}
H_{11}(s) & H_{12}(s)\\H_{21}(s) & H_{22}(s)
    \end{array}
  \right]. 
\end{equation}
Also, 
\begin{equation}
\Phi(s) \triangleq
       \left[
       \begin{array}{cc}
\Psi(s) &  -G_{11}(s)(I_n+\Sigma_u^T)\\
-(I_n+\Sigma_u^T)G_{11}(s)^H  &~~ \Sigma_u^T+2I_n   
       \end{array}       \right],  \label{eq:75} 
\end{equation}
where  
\begin{equation}
  \label{eq:47}
\Psi(s)\triangleq G_{11}(s)\Sigma_u^TG_{11}(s)^H
       +G_{12}(s)\Sigma_w^TG_{12}(s)^H.   
\end{equation}

The problem of optimal coherent equalization posed in~\cite{UJ2} seeks to
infimize  an upper bound on the largest 
eigenvalue of the Hermitian matrix $P_e(i\omega)$ over the entire range of
frequency $\omega$. Formally, it is the optimization problem in which one
wants to find 
\begin{equation}
  \label{eq:6'}
\gamma_\circ\triangleq \inf \gamma  
\end{equation}
over the set of $\gamma>0$ for which
there exists an admissible $H(s)\in \mathcal{H}_p$ such that
\begin{equation}
  \label{eq:6'.sub}
P_e(i\omega)<\gamma^2 I_n \quad \forall\omega.
\end{equation}
In the sequel, we will refer to the constant $\gamma_\circ$ defined
in~(\ref{eq:6'}) and (depending 
on the context) the corresponding matrix $\gamma_\circ^2 I_n$ as an optimal
equalization performance. A $\gamma\ge \gamma_\circ$ and an
admissible $H(s)\in \mathcal{H}_p$ for
which~(\ref{eq:6'.sub}) holds are regarded as a guaranteed cost feasible
solution of problem~(\ref{eq:6'}). 

The objective of the optimization problem~(\ref{eq:6'}) is similar 
to the classical $H_\infty$ filtering objective~\cite{HSK-1999}, except it
employs the largest eigenvalue of $P_e(i\omega)$ instead of a 
singular value of the disturbance-to-error transfer function. As in
$H_\infty$ filtering, the infimum in~(\ref{eq:6'}) may not be attained in
general, therefore of practical interest are feasible solutions with a
bound $\gamma^2$ on the power spectrum density $P_e(i\omega)$ which closely
approximates $\gamma_\circ^2$. Following the analogy with $H_\infty$
filtering, when the gap between $\gamma$ and  $\gamma_\circ$ is acceptably
small, such feasible solutions are termed \emph{suboptimal} or \emph{near
optimal}, since they guarantee an acceptable mean-square mismatch between
the input field in the system and the output of the filter, in the sense of
Definition~\ref{def.error}. The same 
terminology will apply
to other optimization problems considered in the paper.

\section{A convex relaxation and a physically realizable upper
  bound}\label{relaxation} 

Unlike $H_\infty$ filtering, admissible filters for the optimization
problem~(\ref{eq:6'}) are subject to the physical realizability
constraint~(\ref{eq:60}) which is not convex in general.  
As a step towards solving the constrained optimization
problem~(\ref{eq:6'}), reference~\cite{UJ2} introduced the following
auxiliary optimization problem 
\begin{equation}
  \label{eq:7}
 \gamma_\circ'\triangleq \inf\{\gamma>0:
\mathcal{H}_{11,\gamma}\neq\emptyset\},
\end{equation}
where $\mathcal{H}_{11,\gamma}$ is the set of $n\times n_y$
transfer functions  
$H_{11}(s)$ which satisfy~(\ref{eq:6'.sub}) as well as
the following conditions
\begin{enumerate}[(H1)]
\item
All poles of $H_{11}(s)$ are in the open left half-plane of the complex
plane, and $H_{11}(s)$ is analytic in a half-plane
$\mathrm{Re}s>-\tau$ ($\exists \tau>0$); 
\item 
  \begin{equation}
    \label{eq:1.UJ2}
    H_{11}(i\omega)H_{11}(i\omega)^\dagger \le I_n \quad \forall \omega\in
    \mathbf{R}. 
  \end{equation}
\end{enumerate}
This problem can be regarded as a convex relaxation of~(\ref{eq:6'}). 
Indeed, for any admissible $H(s)$ of problem~(\ref{eq:6'}), its (1,1) block
denoted $H_{11}(s)$
satisfies conditions (H1),
(H2). Therefore, the value $\gamma_\circ'$ of problem~(\ref{eq:7}) is a lower bound on
the optimal equalization performance $\gamma_\circ$:
\begin{equation}
  \label{eq:9}
  \gamma_\circ\ge \gamma_\circ'.
\end{equation}

The constant $\gamma_\circ'$ indicates a fundamental performance limitation
of the 
considered coherent equalization scheme, since according 
to~(\ref{eq:9}) no admissible equalizer can approximate $u$ with $\hat u$
so that the corresponding power spectrum density $P_e(i\omega)$ is less
than $(\gamma_\circ')^2I_n$ across the entire range of frequencies.  
It was also shown in~\cite{UJ2} that, given a $\gamma\ge \gamma_\circ'$, if
$H_{11}\in\mathcal{H}_{11,\gamma}$ has the additional property 
\begin{enumerate}[(H3)]
\item
The normal rank\footnote{A non-negative
  integer $r$ is the normal rank of a rational function $X(s)$ if (a)
  $X$ has at least one subminor of order $r$ which does not vanish
  identically,  and (b) all minors of order greater than $r$ vanish
  identically~\cite{Youla-1961}.} of the following matrices does not change
on the finite imaginary axis $i\omega$: 
\begin{eqnarray}
  \label{eq:26}
  Z_1(s)&=&I_n-H_{11}(s)H_{11}(s)^H, \nonumber \\ 
  Z_2(s)&=&I_{n_y}-H_{11}(s)^HH_{11}(s),
\end{eqnarray}
\end{enumerate}
then an admissible physically realizable $(n+r)\times (n+r)$ matrix
transfer function $H(s)$ can be constructed  for which (\ref{eq:6'.sub})
holds with the same $\gamma$
; also see Section~\ref{X}. Here $r$ is the normal rank of 
$Z_2$. 
However, in general, 
the elements of $\mathcal{H}_{11,\gamma}$ are not guaranteed to 
satisfy condition (H3).  
Hence, the existence of an admissible suboptimal 
coherent equalizer cannot be
generally inferred from problem~(\ref{eq:7}).
This leaves a gap between the underlying
problem~(\ref{eq:6'}) and the convex relaxation~(\ref{eq:7}).

To address this gap, we replace the inequality~(\ref{eq:1.UJ2}) with
a strict inequality.
We now show that the relaxation of the
problem~(\ref{eq:6'}) obtained this way yields an \emph{upper} bound on
the optimal equalization performance $\gamma_\circ$, and that suboptimal
solutions of the modified auxiliary 
problem are \emph{guaranteed} to satisfy condition (H3). As a result, a
physically realizable filter which guarantees equalization
accuracy~(\ref{eq:6'.sub}) can be constructed from every suboptimal
solution of the new auxiliary problem, as will be explained in Section~\ref{X}. 

Given $\gamma>0$, define
the set $\mathcal{H}_{11,\gamma}^-$ consisting of proper rational $n\times n_y$ 
transfer functions which satisfy (H1) and~(\ref{eq:6'.sub}) as well as
the following condition
  \begin{equation}
    \label{eq:1}
    H_{11}(i\omega)H_{11}(i\omega)^\dagger < I_n \quad \forall \omega\in
    \bar{\mathbf{R}}.
  \end{equation}
Here $\bar{\mathbf{R}}$ denotes the closed
real axis: $\bar{\mathbf{R}}\triangleq \mathbf{R} \cup\{\pm\infty\}$.  

In view of (\ref{eq:1}), we will call transfer functions of the set
$\mathcal{H}_{11,\gamma}^-$ strictly contractive, or contractive for
short. Indeed, (\ref{eq:1}) is equivalent to 
$\|H_{11}\|_\infty=\sup_\omega\|H_{11}(i\omega)\|<1$. Transfer
functions with this property correspond to causal contractive mappings
$L_2[-\infty,\infty)\to L_2[-\infty,\infty)$. 

\begin{remark}\label{R.contractive}
Mathematically, condition~(\ref{eq:1.UJ2}) describes the closed unit ball in
$H_\infty$, and condition~(\ref{eq:1}) describes its open
interior. Therefore, contractiveness does not preclude $H_{11}(i\omega)$
from being arbitrarily close to the boundary of this unit ball. Physically, when
the channel $y\to \hat u$ of the equalizer is contractive, the filter attenuates
$y(t)$, and along with it, the degrading contribution of the thermal noise
of the system environment is attenuated. 
\hfill$\Box$
\end{remark}

Note that the set $\mathcal{H}_{11,\gamma}^-$ is not empty for any $\gamma>0$
such that $\gamma^2 I\ge 2I+\Sigma_u^T$; it contains the transfer function
$H_{11}(s)\equiv 0$. Thus the following constant is well defined,
\begin{equation}
  \label{eq:6'''}
\gamma_\circ''\triangleq \inf\{\gamma>0\colon
  \mathcal{H}_{11,\gamma}^-\neq\emptyset\}. 
\end{equation}

\begin{theorem}
  \label{P1a>P1}
Suppose $\gamma>0$ is such that $\mathcal{H}_{11,\gamma}^-$ is not
empty. Then every element $H_{11}(s)$ of
$\mathcal{H}_{11,\gamma}^-$ satisfies all three conditions (H1)-(H3). As a
result, an admissible physically realizable filter transfer function
$H(s)\in\mathcal{H}_p$ can be constructed using this $H_{11}(s)$ as the
(1,1) block of $H(s)$ which guarantees the equalization
performance~(\ref{eq:6'.sub}). The infimal equalization 
performance among all admissible $H(s)$ constructed this way
is equal to the value of $\gamma_\circ''$ 
defined in~(\ref{eq:6'''}), therefore it holds that
  \begin{equation}
    \label{eq:2}
   \gamma_\circ''\ge  \gamma_\circ. 
  \end{equation}
\end{theorem}

\emph{Proof: }
Let  $H_{11}(s)\in \mathcal{H}_{11,\gamma}^-$, and so
$\gamma\ge\gamma_\circ''$. Then by definition, 
$H_{11}(s)$ satisfies condition (H1). (H2) follows
from~(\ref{eq:1}) trivially. Furthermore, owing to~(\ref{eq:1}), the matrices 
$Z_1(i\omega)$, $Z_2(i\omega)$ defined in~(\ref{eq:26}) are nonsingular for
every $\omega\in\mathbf{R}$. Therefore, $\det Z_1(s)$, $\det Z_2(s)$ do not
vanish on the finite imaginary axis, and consequently condition (H3) is
also satisfied. This proves the first claim of the theorem. 

Furthermore, since $H_{11}(s)$ has been shown to satisfy all three
conditions (H1)-(H3), 
the result of Theorem~1 in~\cite{UJ2} can now be applied to construct an  
admissible physically realizable filter transfer function
$H(s)\in\mathcal{H}_p$ using this $H_{11}(s)$; Section~\ref{X} will explain
this in detail. This transfer function will
satisfy condition~(\ref{eq:6'.sub}), since by construction the (1,1) block of
$H(s)$ is equal to the selected element $H_{11}(s)$ of the set
$\mathcal{H}_{11,\gamma}^-$ (see Section~\ref{X}), and so  
$P_e(i\omega)<\gamma^2I$ according to the definition of this set. 
Thus, we
conclude that $\gamma\ge \gamma_\circ$.   
Taking infimum over $\gamma$ leads to~(\ref{eq:2}). 
\hfill$\Box$

Note that both problems~(\ref{eq:7}) and (\ref{eq:6'''}) are
convex. In the next sections, we discuss a method to solve these auxiliary
convex problems. However, 
Theorem~\ref{P1a>P1} shows that unlike~(\ref{eq:7}), \emph{any} feasible
solution to problem~(\ref{eq:6'''}) (i.e., a constant
$\gamma>\gamma_\circ''$ 
and a corresponding transfer function $H_{11}\in\mathcal{H}_{11,\gamma}^-$)  
can be used to construct a 
feasible physically realizable transfer function
$H(s)$, and one can ascertain that the power spectrum density of the
corresponding equalization error will be bounded from above by $\gamma^2I_n$.  
That is, the equalizers
constructed from elements of the set $\mathcal{H}_{11,\gamma}^-$ provide a
\emph{bona fide} guarantee of the mean-square accuracy. 
Furthermore, by minimizing over $\gamma$ such a filter can be constructed
to guarantee the mean-square accuracy 
arbitrarily close to $\gamma_\circ''$. In general, since
$\gamma_\circ''\ge\gamma_\circ$, there may remain some performance gap
between the equalization performance of this filter and the theoretically
smallest equalization performance symbolized by $\gamma_\circ$. In the next
section, a condition will be obtained under which the
relaxation~(\ref{eq:6'''}) is exact and this gap vanishes.

\section{Suboptimal solutions to the auxiliary problems}\label{framework}

In this section, we present a method to characterize suboptimal solutions to
the auxiliary optimization problem~(\ref{eq:6'''}) by recasting the latter
problem as a finite dimensional optimization problem.

Problem~(\ref{eq:6'''}) involves frequency domain constraints~(\ref{eq:6'.sub})
and~(\ref{eq:1}). A 
powerful approach to solving optimization problems with this type of
constraints is based on
the Kalman-Yakubovich-Popov (KYP) 
lemma~\cite{DP-2000,Rantzer-1996,Yakubovich-1974}. It replaces the infinite
number of constraints 
indexed by the frequency with constraints expressed in terms of constant
matrices, making the problem finite dimensional. However,  
the standard KYP lemma is not directly applicable to problem
(\ref{eq:6'''}), since it involves two frequency domain 
constraints. 
Therefore, in order
to apply this approach to~(\ref{eq:6'''}), we
adopt a sufficient condition on the system~(\ref{dyn}) under which
the constraint~(\ref{eq:1}) becomes superfluous when
(\ref{eq:6'.sub}) is satisfied. In fact, this condition is a
simple modification of~\cite[Theorem~2]{UJ2}.

Given $\gamma>0$, let the set $\mathcal{H}_{11,\gamma}^\infty$
consist of $n\times n_y$ matrix elements of the Hardy space $H_\infty$
which satisfy condition (H1) and
inequality~(\ref{eq:6'.sub}). Then consider the optimization problem
\begin{equation}
  \label{eq:71.nc}
\gamma_*=\inf \{\gamma>0: \mathcal{H}_{11,\gamma}^\infty\neq\emptyset\}.
\end{equation} 

Unlike problems~(\ref{eq:7}) and~(\ref{eq:6'''}), problem~(\ref{eq:71.nc})
does not impose magnitude constraints on $H_{11}(i\omega)$. Furthermore, it
holds that 
$\mathcal{H}_{11,\gamma}^-\subseteq \mathcal{H}_{11,\gamma} \subseteq \mathcal{H}_{11,\gamma}^\infty$,
therefore the optimal values of all four optimization problems considered
in this paper are in the ordered relation, 
\begin{equation}
  \label{eq:15}
 \gamma_*\le \gamma_\circ' \le \gamma_\circ \le
\gamma_\circ''.  
\end{equation}
The leftmost inequality follows from the aforementioned set inclusion,
while the other two inequalities have been established earlier;
see~(\ref{eq:9}) and~(\ref{eq:2}). The following theorem gives a sufficient condition under which this chain of
inequalities is a chain of identities, and so the 
relaxations~(\ref{eq:6'''}) and~(\ref{eq:7}) become exact.   

\begin{theorem}\label{SDP.primal.LMI}
Suppose $\gamma>0$ is such that
$\mathcal{H}_{11,\gamma}^-\neq\emptyset$. 
If there exists $\theta >0$ such
that  
\begin{equation}
\label{eq:21}
\theta \left(\Phi(i\omega)-\gamma^2
\left[
  \begin{array}{cc}
    0&~0 \\
    0&~ I_n
  \end{array}
\right]\right)-\left[
  \begin{array}{cc}
    I_{n_y}&0 \\
    0& -I_n
  \end{array}
\right]>0 \quad \forall\omega\in\bar{\mathbf{R}},
\end{equation}
then $\mathcal{H}_{11,\gamma}^-= \mathcal{H}_{11,\gamma}^\infty$. 
Consequently, the optimal equalization performance $\gamma_\circ$ defined in~(\ref{eq:6'}) is equal
to the values of all three convex relaxation
problems~(\ref{eq:71.nc}),~(\ref{eq:7}) and~(\ref{eq:6'''}),   
\begin{equation}
  \label{eq:11}
\gamma_\circ = \gamma_*= \gamma_\circ' = \gamma_\circ''.
\end{equation}
\end{theorem}

{\emph{Proof: }
The proof of the set equality
$\mathcal{H}_{11,\gamma}^-= \mathcal{H}_{11,\gamma}^\infty$ is
identical to the proof of Theorem~2 in~\cite{UJ2}, except for the obvious
changes due to the fact that~(\ref{eq:21}) is a strict inequality. It must
hold strictly to guarantee that (\ref{eq:1}) is satisfied 
when $\omega\to\pm\infty$. Also,
since these sets are equal then  $\gamma_*= \gamma_\circ''$. The
identity~(\ref{eq:11}) then follows immediately from~(\ref{eq:15}). 
\hfill$\Box$

\begin{remark}\label{Rem.about.21}
Condition (\ref{eq:21}) relates the gain of the
system~(\ref{dyn}) to the intensities $\Sigma_u$, $\Sigma_w$ of the signal
and noise 
inputs. If this condition holds for a certain environment noise intensity
$\Sigma_w^0$, then it holds for any $\Sigma_w\ge \Sigma_w^0$; this is the
low SNR situation to which we alluded in the introduction. Furthermore, if
(\ref{eq:21}) holds for a 
certain $\gamma^0$, then it holds for any $\gamma\in (0,\gamma^0]$. 
Nevertheless, the corresponding set $\mathcal{H}_{11,\gamma}^-$ may be
empty for some such $\gamma$.   
\hfill$\Box$
\end{remark}

From Theorem~\ref{SDP.primal.LMI}, it follows that under
condition~(\ref{eq:21}), any suboptimal transfer function $H_{11}(s)$ of
problem~(\ref{eq:71.nc}) is also suboptimal for the auxiliary 
problem~(\ref{eq:6'''}). Furthermore, recalling the discussion in the last
paragraph of Section~\ref{relaxation}, equation~(\ref{eq:11}) shows that
any near optimal transfer function $H_{11}(s)$ of problem~(\ref{eq:71.nc})
can also be used to construct a near optimal admissible $H(s)$ for the
underlying equalization problem~(\ref{eq:6'}). Thus if
condition~(\ref{eq:21}) holds, constraint~(\ref{eq:1}) can be dropped while
searching for a solution to problem~(\ref{eq:6'''}). This observation
reduces computing a near optimal solution to the auxiliary
problem~(\ref{eq:6'''}) to 
computing a constant $\gamma$ approximating the value $\gamma_*$ of
problem~(\ref{eq:71.nc}) with a desired accuracy
and a corresponding transfer function $H_{11}\in \mathcal{H}_{11,\gamma}^\infty$.
We now present a solution to the latter problem. It
relies on the following technical assumption.

\begin{assumption}
  \label{A.rho}
There exists a constant $\lambda\ge 0$ such that the $(n_y+n)\times
(n_y+n)$ rational matrix transfer function 
\begin{equation}
  \label{eq:39}
  \Phi_\lambda(s)=\Phi(s)+
  \left[
    \begin{array}{cc}
      0 & 0 \\
      0 &~~  \lambda^2 I_n
    \end{array}
  \right]
\end{equation}
admits a spectral factorization
 \begin{equation}
  \label{eq:68}
  \Phi_\lambda(s)=\Upsilon_\lambda(s)\Upsilon_\lambda(s)^H, 
\end{equation}
where a $(n_y+n)\times p$ rational transfer matrix $\Upsilon_\lambda(s)$ has all
its poles in the left half-plane $\mathrm{Re}s<-\tau$ and is analytic in
$\mathrm{Re}s>-\tau$ ($\exists \tau>0$).
\end{assumption}

The purpose of adding the second term on the right-hand side
of~(\ref{eq:39}) is to make $\Phi_\lambda$ positive
semidefinite on the imaginary axis. In this case, a spectral factorization of
$\Phi_\lambda(s)$ exists~\cite{Anderson-1967,CALM-1997,Youla-1961}. An algebraic
method to compute spectral factors can also 
be found in~\cite{Anderson-1967}. 
The following lemma gives a simple sufficient condition for
$\Phi_\lambda(i\omega)\ge 0$ for all $\omega$.  

\begin{lemma}\label{L.Phi>=0}
  If there exists $\lambda>0$ such that 
  \begin{equation}
    \label{eq:20}
    \Sigma_u^T\ge
    (\Sigma_u^T+I_n)(\Sigma_u^T+(2+\lambda^2)I_n)^{-1}(\Sigma_u^T+I_n), 
  \end{equation}
then $\Phi_\lambda(i\omega)\ge 0$ for all $\omega$.
\end{lemma}

\emph{Proof: }
Since $\Sigma_w^T\ge 0$, it follows from~(\ref{eq:20}) that
\begin{eqnarray*}
\Psi(i\omega)\ge
G_{11}(i\omega)(\Sigma_u^T+I_n)(\Sigma_u^T+(2+\lambda^2)I_n)^{-1}
  \\ 
\times (\Sigma_u^T+I_n)G_{11}(i\omega)^\dagger.
\end{eqnarray*}
The claim of the lemma then follows using the Schur complement.
\hfill$\Box$ 

\begin{remark}
It is easy to show that~(\ref{eq:20}) holds with some sufficiently large
$\lambda$ when $\Sigma_u^T> 0$. Indeed if $\Sigma_u^T> 0$, then
$\Sigma_u^T+I_n>I_n$. This implies that
$(\Sigma_u^T+I_n)^{-1}>(\Sigma_u^T+I_n)^{-2}$. Consequently it 
is possible to choose a large enough $\lambda^2$ so that
\[
(\Sigma_u^T+I_n)^{-1}-(\Sigma_u^T+(2+\lambda^2)I_n)^{-1}>(\Sigma_u^T+I_n)^{-2},
\]
which implies~(\ref{eq:20}).
\hfill$\Box$  
\end{remark}   

\begin{lemma}\label{L.gam=gambar}
  Consider the optimization problem
\begin{eqnarray}
  \label{eq:6}
&&  \bar\gamma_*\triangleq \inf \bar \gamma \\
&&\mbox{subject to $\bar\gamma>\lambda$ and} \nonumber \\
  \label{eq:3}
&&\left[
  \begin{array}{cc}
    H_{11}(i\omega) & I_n
  \end{array}
  \right]\Upsilon_\lambda(i\omega)\Upsilon_\lambda(i\omega)^\dagger  \left[
\begin{array}{c}
    H_{11}(i\omega)^\dagger \\ I_n
  \end{array}
  \right]<\bar\gamma^2I_n \nonumber \\
&& \phantom{\left[
  \begin{array}{cc}
    H_{11}(i\omega) & I_n
  \end{array}
  \right]\Upsilon_\lambda(i\omega)\Upsilon_\lambda(i\omega)^\dagger  \left[
\begin{array}{c}
    H_{11}(i\omega)^\dagger
  \end{array}
  \right]}\forall \omega\in\mathbf{R},
\end{eqnarray}
where $H_{11}$ satisfies condition (H1). Then
\begin{equation}
  \label{eq:62}
  \bar\gamma_*^2=\gamma_*^2+\lambda^2.
\end{equation}
Furthermore, for any $\gamma>\gamma_*$, $H_{11}\in H_{11,\gamma}^\infty$
if and only if~(\ref{eq:3}) holds with
$\bar\gamma=(\gamma^2+\lambda^2)^{1/2}>\bar\gamma_*$. 
\end{lemma}

\emph{Proof: }
For any $\gamma>\gamma_*$, $\mathcal{H}_{11,\gamma}^\infty\neq
\emptyset$. Then if 
$H_{11}\in\mathcal{H}_{11,\gamma}^\infty$, then (\ref{eq:6'.sub})
implies~(\ref{eq:3}) with $\bar\gamma=(\gamma^2+\lambda^2)^{1/2}>\lambda$. Hence,
such $H_{11}$ and $\bar\gamma$ are in the feasible set of
problem~(\ref{eq:6}), and $\bar\gamma\ge \bar\gamma_*$, i.e., 
$\gamma^2+\lambda^2\ge \bar\gamma_*^2$. This implies that
$\gamma_*^2+\lambda^2\ge \bar\gamma_*^2$. 

Conversely, 
consider $\bar \gamma>\bar \gamma_*$,  and let $H_{11}$ satisfy (H1)
and~(\ref{eq:3}). Since by definition $\bar\gamma_*\ge \lambda$, then 
$\gamma=(\bar\gamma^2-\lambda^2)^{1/2}>0$ is well defined. Then it follows
from~(\ref{eq:3}) that  $H_{11}$ satisfies (\ref{eq:6'.sub}). Thus we
conclude that $H_{11}\in\mathcal{H}_{11,\gamma}^\infty$, i.e., 
$\mathcal{H}_{11,\gamma}^\infty\neq \emptyset$. This implies that
$\gamma\ge \gamma_*$, hence $\bar\gamma^2\ge \gamma_*^2+\lambda^2$. Taking
infimum over $\bar\gamma$ yields $\bar\gamma_*^2\ge \gamma_*^2+\lambda^2$.
Combining this inequality with the previously established inequality
$\gamma_*^2+\lambda^2\ge \bar\gamma_*^2$, we conclude that (\ref{eq:62})
holds true.  
\hfill$\Box$

Lemma~\ref{L.gam=gambar} shows that problems~(\ref{eq:71.nc})
and~(\ref{eq:6}) are 
equivalent. However, problem~(\ref{eq:6}) has an advantage in that the
factorization~(\ref{eq:39}) facilitates a solution in the tractable form of 
a semidefinite program, as we now demonstrate. 

Let $m$ be the McMillan degree of $\Upsilon_\lambda(s)$. Introduce a
minimal state-space representations of this transfer function, 
  \begin{equation}
    \label{eq:17}
    \Upsilon_\lambda(s)\sim
    \left[
    \begin{array}{c|c}
     A_\lambda & B_\lambda \\
\hline
     C_{1,\lambda} & D_{1,\lambda} \\ 
     C_{2,\lambda} & D_{2,\lambda} \\ 
    \end{array}
    \right].
  \end{equation}
  The complex matrices $A_\lambda$, $B_\lambda$, $C_{1,\lambda}$, $C_{2,\lambda}$,
  $D_{1,\lambda}$, $D_{2,\lambda}$ have dimensions $m\times m$, $m\times p$,
  $n_y\times m$, $n\times m$, $n_y\times p$, $n\times p$, respectively; they are
  obtained from the factorization~(\ref{eq:68}) and are determined by the
  characteristics of the channel transfer function  $G(s)$ and matrices
  $\Sigma_u$, $\Sigma_w$ of the inputs $u$, $w$. They represent the
  data of the filter design problem. By Assumption~\ref{A.rho},
  $A_\lambda$ is a Hurwitz matrix.

Also introduce a state-space
  representation for the transfer function $H_{11}(s)$, 
  \begin{equation}
\label{eq:45}
     H_{11}(s)\sim
     \left[
     \begin{array}{c|c}
      A_{11} & B_{11} \\
 \hline
      C_{11} & J_{11} 
     \end{array}
     \right].
   \end{equation}
The complex matrices $A_{11}$, $B_{11}$, $C_{11}$, $J_{11}$
will be the decision variables of the semidefinite program in
question, along with $\bar\gamma$. They have dimensions 
$m_{11}\times  
m_{11}$, $m_{11}\times n_y$, $n\times m_{11}$, $n\times n_y$. The dimension of the state-space representation $m_{11}$ is assumed to be fixed; otherwise it is chosen
arbitrarily. Also, since $H_{11}(s)$ is sought in the class of stable transfer
functions, we will require that $A_{11}$ must be Hurwitz. 

For our derivation, it is convenient to introduce
\begin{eqnarray}
  \label{eq:4}
  \bar H_{11}(s)=H_{11}(s^*)^\dagger, \quad 
  \bar \Upsilon_\lambda(s)=\Upsilon_\lambda(s^*)^\dagger. 
\end{eqnarray}
The transfer functions $\bar H_{11}(s)$,
$\bar\Upsilon_\lambda(s)$ inherit the stability and analiticity 
properties of $H_{11}(s)$, $\Upsilon_\lambda(s)$. Their state-space
realizations can be obtained from the state-space realizations of
$H_{11}(s)$, $\Upsilon_\lambda(s)$, as follows 
  \begin{equation}
    \label{eq:17bar}
    \bar\Upsilon_\lambda(s)\sim
    \left[
    \begin{array}{c|cc}
     A_\lambda^\dagger & C_{1,\lambda}^\dagger & C_{2,\lambda}^\dagger \\
\hline
     B_{\lambda}^\dagger & D_{1,\lambda}^\dagger & D_{2,\lambda}^\dagger 
    \end{array}
    \right],
\quad 
    \bar H_{11}(s)\sim
    \left[
    \begin{array}{c|c}
     A_{11}^\dagger & C_{11}^\dagger \\
\hline
     B_{11}^\dagger & J_{11}^\dagger 
    \end{array}
    \right].
  \end{equation}
Using this notation, (\ref{eq:3}) is expressed as
\begin{eqnarray}
  \label{eq:3a}
&&  \left[
  \begin{array}{cc}
    \bar H_{11}(-i\omega)^\dagger & I_n
  \end{array}
  \right]\bar\Upsilon_\lambda(-i\omega)^\dagger\bar\Upsilon_\lambda(-i\omega)  \left[
\begin{array}{c}
    \bar H_{11}(-i\omega) \\ I_n
  \end{array}
  \right]<\bar\gamma^2 I_n \nonumber \\
&& \hspace{5cm}\forall \omega\in\mathbf{R}. \quad
\end{eqnarray}
Since $\omega$ is arbitrary, one can replace $i\omega$ with $-i\omega$
in~(\ref{eq:3a}), to express (\ref{eq:3a}) in a 
form routinely arising in problems of $H_\infty$ control and filtering:
\begin{eqnarray}
\label{eq:8}
&& T_\lambda(i\omega)^\dagger T_\lambda(i\omega)
  <\bar\gamma^2I_n \quad
 \forall \omega\in\mathbf{R}, \\
&& T_\lambda(s)\triangleq \bar \Upsilon_\lambda(s)\left[
\begin{array}{c}
   \bar  H_{11}(s) \\ I_n
  \end{array}
  \right]. \nonumber 
\end{eqnarray}
Inequality (\ref{eq:8}) can be seen as a bound on singular values of the
frequency response of the transfer function $T_\lambda(s)$. This makes it
possible to apply results from the $H_\infty$ control theory~\cite{DP-2000} to 
equivalently express~(\ref{eq:8}) --- and consequently~(\ref{eq:3}) --- in
the form of tractable matrix inequalities\footnote{It worth pointing out
  that the $H_\infty$ problem in hand is a somewhat less common problem of
  computing a prefilter.}.

The first step is to combine
the state-space representations~(\ref{eq:17bar}) into the state-space model
of $T_\lambda(s)$:
  \begin{eqnarray}
    \label{eq:33}
T_\lambda(s)&\sim& 
    \left[
    \begin{array}{c|c}
     \hat A~ &~ \hat B \\
\hline
     \hat C~ &~ \hat D
    \end{array}
              \right] \nonumber \\
&=&\left[
    \begin{array}{cc|c}
      A_\lambda^\dagger & C_{1,\lambda}^\dagger B_{11}^\dagger & 
                C_{1,\lambda}^\dagger J_{11}^\dagger+C_{2,\lambda}^\dagger \\
      0 & A_{11}^\dagger & C_{11}^\dagger 
      \\
\hline
    B_{\lambda}^\dagger & ~~ D_{1,\lambda}^\dagger B_{11}^\dagger~ &~  
                   D_{1,\lambda}^\dagger J_{11}^\dagger+D_{2,\lambda}^\dagger \\
    \end{array}
    \right].
  \end{eqnarray}
Since $H_{11}$, $\Upsilon_\lambda$ are stable, the matrices $A_\lambda^\dagger$,
$A_{11}^\dagger$ and $\hat A= \left[
    \begin{array}{cc}
      A_\lambda^\dagger & C_{1,\lambda}^\dagger B_{11}^\dagger \\
      0 & A_{11}^\dagger 
    \end{array}
    \right]$ 
are stable. Then using the KYP
lemma~\cite{DP-2000,Rantzer-1996,Yakubovich-1974} and the 
Schur complement~\cite{HZ-2005}, we obtain that inequality~(\ref{eq:8}) (and
consequently~(\ref{eq:3})) holds if and only if there exists a
Hermitian positive definite
matrix $\hat X$ such that
\begin{equation}
      \label{eq:35}
      \left[
        \begin{array}{ccc}
       \hat A^\dagger \hat X + \hat X \hat A~  &~ \hat X\hat B  & \hat C^\dagger
          \\
\hat B^\dagger\hat X & -\bar\gamma^2I_n & \hat D^\dagger \\
\hat C &  \hat D &  -I_p
        \end{array}
      \right]<0;
\end{equation}
e.g., see~\cite[Corollary~7.5]{DP-2000}.

Next, define the matrices 
\begin{eqnarray}
  \label{eq:25}
&&  K=
  \left[
  \begin{array}{cc}
    A_{11}^\dagger & C_{11}^\dagger \\
   B_{11}^\dagger & J_{11}^\dagger
  \end{array}
  \right], \quad
\bar A=  \left[
  \begin{array}{cc}
    A_\lambda^\dagger &~ 0_{m\times m_{11}} \\
   0_{m_{11}\times m} &~ 0_{m_{11}\times m_{11}}
  \end{array}
  \right], \nonumber\\
&& \bar B= \left[
  \begin{array}{c}
    C_{2,\lambda}^\dagger  \\
   0_{m_{11}\times n} 
  \end{array}
  \right], \quad
\underline{B}=  \left[
  \begin{array}{cc}
   0_{m\times m_{11}} &~ C_{1,\lambda}^\dagger  \\
   I_{m_{11}} & 0_{m_{11}\times n}
  \end{array}
  \right], \nonumber \\
&& \underline{C}=  \left[
  \begin{array}{cc}
   0_{m_{11}\times m} &~ I_{m_{11}}  \\
   0_{n\times m} & 0_{n\times m_{11}}
  \end{array}
  \right], \quad
\underline{D}_{12}=  \left[
  \begin{array}{cc}
   0_{p\times m_{11}} &~ D_{1,\lambda}^\dagger  \\
  \end{array}
  \right], \nonumber \\
&&
\bar C= \left[
    \begin{array}{cc}
    B_\lambda^\dagger & 0_{p\times m_{11}} \\
    \end{array}
    \right], \quad \bar D=D_{2,\lambda}^\dagger, \nonumber \\
&& \underline{D}_{21}=  \left[
  \begin{array}{c}
   0_{m_{11}\times n}  \\
   I_n 
  \end{array}
  \right].
\end{eqnarray}
This notation reveals that the state-space matrices of the transfer
function $T_\lambda$ are affine in $K$:
\begin{eqnarray}
  \label{eq:31}
  \hat A&=& \bar A+\underline{B}K\underline{C}, \quad 
  \hat B= \bar B+\underline{B}K\underline{D}_{21}, \nonumber \\
  \hat C&=& \bar C+\underline{D}_{12}K\underline{C}, \quad 
  \hat D= \bar D+\underline{D}_{12}K\underline{D}_{21}.
\end{eqnarray}
Using this notation, (\ref{eq:35}) can be written as
\begin{equation}
  \label{eq:40}
  \Sigma_{\hat X}+\Pi^\dagger K^\dagger \Lambda_{\hat X} +\Lambda_{\hat X}^\dagger K \Pi<0,  
\end{equation}
where
\begin{eqnarray*}
  \Sigma_{\hat X}&=&\left[
        \begin{array}{ccc}
          \bar A^\dagger \hat X + \hat X \bar A~  &~ \hat X\bar B  & \bar C^\dagger
          \\
\bar B^\dagger\hat X & -\bar\gamma^2I_n & \bar D^\dagger \\
\bar C &  \bar D &  -I_{p}
        \end{array}
      \right], \\
 \Lambda_{\hat X}&=&\left[
        \begin{array}{ccc}
          \underline{B}^\dagger \hat X  &~ 0_{(m_{11}+n_y)\times n}  &~ \underline{D}_{12}^\dagger
        \end{array}
      \right], \\
\Pi&=&\left[
        \begin{array}{ccc}
          \underline{C} &~ \underline{D}_{21} & 0_{(m_{11}+n)\times p}~
        \end{array}
      \right] \\
&=&
\left[
\begin{array}{cccc}
0_{(m_{11}+n)\times m} & I_{m_{11}+n}& 0_{(m_{11}+n)\times p} \\
\end{array}
\right].
\end{eqnarray*}
These observations combined with Theorem~7.10
from~\cite{DP-2000}\footnote{Even though 
  this result is presented in~\cite{DP-2000} for real matrices, the
  derivations carry over to the complex case.} give necessary and
sufficient conditions for the existence of a matrix $K$ 
which satisfies~(\ref{eq:40}) and equivalently, (\ref{eq:35}).
     
\begin{lemma}
  \label{T.LMI}
Suppose Assumption~\ref{A.rho} is satisfied. A constant $\bar\gamma$ and
matrices $A_{11},B_{11},C_{11},J_{11}$ and $\hat X=\hat X^\dagger>0$ for
which~(\ref{eq:40}) holds exist if and only 
if the following matrix inequalities are feasible in the variables
$\bar\gamma^2$, $X_1=X_1^\dagger\in\mathbf{C}^{m\times m}$,
$Y_1=Y_1^\dagger\in\mathbf{C}^{m\times m}$: 
\begin{eqnarray}
&&\bar\gamma^2>\lambda^2,
\label{eq:52}\\
\label{eq:36}
&& X_1>0, \quad Y_1>0,  \\
  \label{eq:55}
&& A_\lambda X_1+X_1A_\lambda^\dagger + B_\lambda B_\lambda^\dagger<0, \\
&&\left[
    \begin{array}{cc}
 N_c^\dagger & 0 \\
 0 & I_n
    \end{array}
  \right] 
  \left[
    \begin{array}{cccc}
A_\lambda^\dagger Y_1 +Y_1 A_\lambda &~ Y_1B_\lambda  & C_{2,\lambda}^\dagger \\
B_\lambda^\dagger Y_1 & -I_p & D_{2,\lambda}^\dagger \\
C_{2,\lambda} &~ D_{2,\lambda} & -\bar\gamma^2 I_n 
    \end{array}
  \right]   \left[
    \begin{array}{cc}
 N_c & 0 \\
 0 & I_n \\
    \end{array}
  \right]<0, \nonumber \\
\label{eq:82} \\
  \label{eq:86}
 && \left[
  \begin{array}{cc}
    X_1 & I \\ I & Y_1
  \end{array}
  \right]\ge 0, \\
 && \mathrm{rank}\left[
  \begin{array}{cc}
    X_1 & I \\ I & Y_1
  \end{array}
  \right]\le m+m_{11},
\label{eq:19} 
\end{eqnarray}
where $N_c$ is a full rank matrix which spans $\mathrm{Ker}
\left[
  \begin{array}{ccc}
    C_{1,\lambda}& D_{1,\lambda}  
  \end{array}
\right]$:
\begin{equation}
\label{eq:76}
\mathrm{Im}N_c=\mathrm{Ker}
\left[
  \begin{array}{ccc}
    C_{1,\lambda}& D_{1,\lambda} 
  \end{array}
\right].
\end{equation}
\end{lemma}

Lemma~\ref{T.LMI} will be used in Section~\ref{X} to compute a
physically realizable $H(s)$. Of course, we are interested in a smallest
$\bar\gamma^2$ for which 
a feasible solution to the equalizer design problem exists. According to
Lemma~\ref{T.LMI}, such a $\bar\gamma^2$ can be obtained as a (near) optimal
solution of the optimization problem  
\begin{equation}
   \label{eq:38}
   \inf\bar\gamma^2 \mbox{ subject
     to~(\ref{eq:52})--(\ref{eq:19}).} 
\end{equation}
From the foregoing discussion, the value of this problem is exactly
equal to $\bar\gamma_*^2$, where $\bar\gamma_*$ is the value  of problem~(\ref{eq:6}). This leads to the following theorem.

\begin{theorem}\label{main.T}
  Suppose Assumption~\ref{A.rho} is satisfied. Let $\bar\gamma$ lie in the
  feasible set of the optimization problem~(\ref{eq:38}). 
If condition~(\ref{eq:21}) is satisfied with
$\gamma=(\bar\gamma^2-\lambda^2)^{1/2}$, then the set
$\mathcal{H}_{11,\gamma}^-$  corresponding to this $\gamma$ is not empty
and $\gamma\ge \gamma_\circ''=\gamma_\circ$. One transfer function
$H_{11}(s)$ which belongs to this set is given by  
\begin{equation}
  \label{eq:51}
  H_{11}(s)=C_{11}(sI-A_{11})^{-1}B_{11}+J_{11},
\end{equation} 
where $A_{11},B_{11},C_{11},J_{11}$ are matrices which render the
inequality~(\ref{eq:40}) true with the chosen $\bar\gamma$ and some $\hat X=\hat
X^\dagger>0$. That is, $\gamma=(\bar\gamma^2-\lambda^2)^{1/2}$ and the transfer
  function $H_{11}(s)$ in~(\ref{eq:51}) lie in the feasible set
of the optimization problem~(\ref{eq:6'''}). 
\end{theorem}

\emph{Proof: }
From Lemma~\ref{T.LMI}, choosing $\bar\gamma$ from the
  feasible set of the optimization problem~(\ref{eq:38}) defined by the
  matrix inequalities (\ref{eq:52})--(\ref{eq:19}) ensures the existence of
  matrices $A_{11}$, $B_{11}$, $C_{11}$,
$J_{11}$ and $\hat X=\hat X^\dagger>0$ which render the
inequality~(\ref{eq:40}) true.  
Since~(\ref{eq:40}) is equivalent to~(\ref{eq:35}), it follows
from~(\ref{eq:35}) that 
\[
\hat A^\dagger \hat X+\hat X\hat A <0.
\]  
Furthermore, since $\hat X$ is positive definite, this inequality implies that
the matrix $\hat A$ is Hurwitz. Consequently, $A_{11}$ is a Hurwitz
matrix since $\hat A$ is upper block-triangular. 
As a result, the transfer function~(\ref{eq:51}) 
satisfies condition (H1). Also,~(\ref{eq:35}) is
equivalent to~(\ref{eq:8}) via the KYP lemma, and these matrix
inequalities are equivalent to~(\ref{eq:6'.sub}), as
was explained in the discussion preceding Lemma~\ref{T.LMI}. This
ensures that~(\ref{eq:6'.sub}) is satisfied with the constructed matrix
$H_{11}(s)$ and $\gamma=(\bar\gamma^2-\lambda^2)^{1/2}$. Therefore, it
follows from Lemma~\ref{L.gam=gambar} that  $H_{11}\in
\mathcal{H}_{11,\gamma}^\infty$.

Next, observe that if $\gamma=(\bar\gamma^2-\lambda^2)^{1/2}$ satisfies the
condition~(\ref{eq:21}) of Theorem~\ref{SDP.primal.LMI}, then according to
this theorem,
$\mathcal{H}_{11,\gamma}^\infty=\mathcal{H}_{11,\gamma}^-$ and
$\gamma_\circ''=\gamma_\circ$. Thus, $H_{11}\in\mathcal{H}_{11,\gamma}^-$
and the latter set is not empty. By definition 
of this set, this implies that $\gamma=(\bar\gamma^2-\lambda^2)^{1/2}\ge
\gamma_\circ''=\gamma_\circ$. This proves the statement of the theorem.    
\hfill$\Box$

Note that while inequalities~(\ref{eq:36})--(\ref{eq:86}) are
linear matrix inequalities in the variables $\bar\gamma^2$, $X_1$, $Y_1$,
the rank constraint~(\ref{eq:19}) is not an LMI. Although solving 
problems with rank constraints has been known to be a difficult problem, 
recently numerical algorithms have been developed for this problem which
perform well in a number of instances~\cite{BCP-2021}. Meanwhile, one can
avoid dealing with intricacies of rank constrained optimization altogether
by letting $m_{11}=m$. Indeed when $m_{11}=m$, (\ref{eq:19}) holds trivially, and the
problem~(\ref{eq:38}) reduces to a semidefinite program, 
\begin{equation}
  \label{eq:38conv}
  \inf\bar\gamma^2 \mbox{ subject to~(\ref{eq:52})--(\ref{eq:86}).}
\end{equation}
This problem is convex, it can be solved efficiently using the existing
numerical algorithms~\cite{LMI,Nesterov-Nemirovskii}.

\begin{corollary}\label{main.T.conv}
Suppose that $m_{11}=m$ and that the conditions of Theorem~\ref{main.T} are
satisfied with a $\bar\gamma$ which belongs to the feasible set of 
the optimization problem~(\ref{eq:38conv}). Then the statement of
Theorem~\ref{main.T} holds true.
\end{corollary}

  The requirement for $\gamma=(\bar\gamma^2-\lambda^2)^{1/2}$
  to satisfy condition~(\ref{eq:21})  of Theorem~\ref{SDP.primal.LMI} 
  ensures that the transfer function $H_{11}(s)$ 
  obtained from problems~(\ref{eq:38}) or~(\ref{eq:38conv})   
  is contractive. Of course, once a feasible $\bar\gamma$
  and the 
  corresponding $H_{11}(s)$ are found from these problems, one can check
  directly whether this $H_{11}(s)$ is contractive. However,
  the frequency domain condition~(\ref{eq:21}) provides an insight into a
  suitable range of $\gamma$ and $\bar\gamma$ \emph{before} one 
  obtains $H_{11}(s)$; see Remark~\ref{Rem.about.21}. 
On the other hand, condition~(\ref{eq:21}) is only
sufficient for 
$H_{11,\gamma}^\infty=H_{11,\gamma}^-$. A contractive transfer function $H_{11}(s)$
may exist even when~(\ref{eq:21})
fails to hold. As we will demonstrate later, in this situation, it may still
be possible to construct a physically realizable equalizer from such
$H_{11}(s)$.

\section{Synthesis of suboptimal physically realizable equalizers}\label{X}     

A procedure for constructing a physically realizable transfer function
$H(s)$ from a suboptimal solution to the auxiliary problem~(\ref{eq:6'''})
was presented in~\cite{UJ2}. 
It was shown  that if a transfer function $H_{11}(s)$
satisfies conditions (H1)--(H3), then the remaining transfer function
matrices $H_{12}(s)$, $H_{21}(s)$, $H_{22}(s)$ of the
partition~(\ref{eq:98a}) which make up an admissible $H(s)$ can
be obtained as follows:
\begin{enumerate}[(i)]
\item
$H_{12}$ is taken to be a left spectral factor of $Z_1(s)$ in~(\ref{eq:26}):
\begin{equation}
Z_1(s)=
H_{12}(s)H_{12}(s)^H.
\label{eq:22}
\end{equation}
\item
Then $H_{21}$ and $H_{22}$ are obtained as
\begin{eqnarray}
  \label{eq:81}
&&H_{21}(s)=U(s)\tilde H_{21}(s), \nonumber \\
&&H_{22}(s)=-U(s)(\tilde H_{21}^{-1}(s))^H H_{11}(s)^H
   H_{12}(s), \quad
\end{eqnarray}
where 
\begin{itemize}
\item
$\tilde H_{21}(s)$ is a right spectral factor of $Z_2(s)$
in~(\ref{eq:26}),  
\begin{equation}
Z_2(s)=
\tilde H_{21}(s)^H  \tilde H_{21}(s); 
\label{eq:85}
\end{equation}
\item
$\tilde H_{21}^{-1}(s)$ is the right inverse of $\tilde
H_{21}(s)$ i.e., an analytic in a right-half-plane $\mathrm{Re}
s>-\tau$ ($\exists\tau>0$) transfer function such that $\tilde H_{21}(s)
\tilde H_{21}^{-1}(s)=I_r$, where $r$ is the normal rank of
$Z_2(s)$~\cite{Youla-1961}; and 
\item
$U(s)$ is a stable, analytic in the closed
right half-plane, paraunitary
$r\times r$ transfer function matrix, chosen to cancel unstable poles of  
$(\tilde H_{21}^{-1}(s))^H H_{11}(s)^H$~\cite{Shaked-1990}.
\end{itemize}
\end{enumerate}

We now connect this procedure with the results of the previous section.   
We have shown in Theorem~\ref{P1a>P1} that any near optimal solution
to~(\ref{eq:6'''}) defines a contractive $H_{11}(s)$ which is
suitable for using in the above procedure, and have also developed 
optimization problems~(\ref{eq:38}) and~(\ref{eq:38conv}) 
which characterize such near optimal solutions in terms of matrix
inequalities; see Theorem~\ref{main.T} and Corollary~\ref{main.T.conv}
respectively. Together, these results lead to 
an algorithm for computing a physically realizable guaranteed cost filter
$H(s)$.

\begin{algo}\label{Alg1}\mbox{}
\begin{enumerate}[1.]
\item 
Select $m_{11}\le m$ and obtain a near optimal
feasible solution (i.e, a constant $\bar\gamma$ and matrices $X_1$, $Y_1$)  
to the optimization problem~(\ref{eq:38}) which approximates the value of
the problem with a suitable accuracy. If $m_{11}=m$ is selected, such
solution can be found 
by solving the convex optimization problem~(\ref{eq:38conv}) to a desired
numerical accuracy.   
These constant $\bar\gamma$ and the matrices $X_1$, $Y_1$ will be used in
the next steps of the procedure.  

\item
Using the found matrices $X_1$, $Y_1$, construct a $m\times m_{11}$ matrix
$X_2$ such that $X_2X_2^\dagger= X_1-Y_1^{-1}$. Such $X_2$ exists since
$X_1-Y_1^{-1}\ge 0$ according to~(\ref{eq:86}). Then construct the
$(m+m_{11})\times (m+m_{11})$ matrix $\hat X$, 
\begin{equation}
  \label{eq:48}
  \hat X=
  \left[
    \begin{array}{cc}
      X_1 & X_2 \\
      X_2^\dagger & I_{m_{11}}
    \end{array}
  \right].
\end{equation}
Note that $\hat X>0$, since by definition of $X_2$, $X_1-X_2X_2^\dagger=Y_1^{-1}
>0$; see~(\ref{eq:36}).
\item
With this $\hat X$ and $\bar\gamma$ found in Step 1, solve the linear
matrix inequality~(\ref{eq:40}) 
for $K$ and
construct the transfer function $H_{11}(s)$ following equation~(\ref{eq:51}).
\item
Using the found $H_{11}(s)$ obtain the transfer functions  $Z_1(s)$ and
$Z_2(s)$ in~(\ref{eq:26}), then compute their spectral factors $H_{12}(s)$
and $\tilde H_{21}(s)$ as per~(\ref{eq:22}),~(\ref{eq:85}). According
to~\cite{CALM-1997}, under condition~(\ref{eq:1}), one suitable choice 
is
\begin{eqnarray}
  \label{eq:27}
  H_{12}(s)&=&-Z_1(i\infty)^{1/2}
      -C_{11}(sI-A_{11})^{-1}L_1,\\
  \tilde H_{21}(s)&=&Z_2(i\infty)^{1/2}
      -L_2(sI-A_{11})^{-1}B_{11}.
\label{eq:34}
\end{eqnarray}
Here, 
\begin{eqnarray*}
&& Z_1(i\infty)=I_n-J_{11}J_{11}^\dagger , \quad
Z_2(i\infty)=I_{n_y}-J_{11}^\dagger J_{11}; \\ 
&& L_1= -(Q_{12}C_{11}^\dagger+B_{11}J_{11}^\dagger)Z_1(i\infty)^{-1/2}, \nonumber \\
&& L_2=Z_2(i\infty)^{-1/2}(Q_{21}B_{11}+C_{11}^\dagger J_{11})^\dagger, \nonumber 
\end{eqnarray*}  
and $Q_{12}$, $Q_{21}$ are Hermitian solutions of the algebraic Riccati equations
\begin{eqnarray}
  \label{eq:23}
\lefteqn{A_{11} Q_{12}+Q_{12}A_{11}^\dagger+(Q_{12}C_{11}^\dagger+B_{11}J_{11}^\dagger)} && \nonumber \\
&&  \times
   Z_1(i\infty)^{-1}(Q_{12}C_{11}^\dagger+B_{11}J_{11}^\dagger)^\dagger
   +B_{11} B_{11}^\dagger= 0, \\
&& \nonumber \\
\lefteqn{Q_{21}A_{11}+A_{11}^\dagger Q_{21}+(Q_{21}B_{11}+C_{11}^\dagger J_{11})} && \nonumber \\
&&  \times Z_2(i\infty)^{-1}(Q_{21}B_{11}+C_{11}^\dagger J_{11})^\dagger+C_{11}^\dagger C_{11} = 0.
\label{eq:24}
\end{eqnarray}
\item
Obtain the remaining transfer functions $H_{21}(s)$ and $H_{22}(s)$
using~(\ref{eq:81}).  
\end{enumerate}
\end{algo}

Steps 1-3 of Algorithm~\ref{Alg1} mirror the procedure for computing a
solution to an $H_\infty$ control problem~\cite{DP-2000}. Here, we apply this
procedure to compute an approximately optimal
bound~(\ref{eq:6'.sub}) on the error power spectrum density and also compute,
via~(\ref{eq:51}),  a
transfer function $H_{11}(s)$ that guarantees this bound. 
Theorem~\ref{main.T} and 
Corollary~\ref{main.T.conv} ensure that the found $H_{11}(s)$ lies in the
feasible set $\mathcal{H}_{11,\gamma}^-$, hence it is suitable for using in
the next 
steps of the algorithm. In particular, the contractiveness property of
$H_{11}$ 
ensures that the matrices $Z_1(i\infty)$ and $Z_2(i\infty)$ are nonsingular
positive definite matrices. Furthermore since the matrix $A_{11}$ produced
at Step 3 is Hurwitz (see the proof of Theorem~\ref{main.T}) and
$\|H_{11}(s)\|_\infty<1$ according to~(\ref{eq:1}), the complex Strict
Bounded Real Lemma~\cite[Theorem~6.2]{MP-2011} ensures that the 
Riccati equations~(\ref{eq:23}) and~(\ref{eq:24}) have Hermitian
solutions. In particular, equation~(\ref{eq:24}) has a nonnegative
definite stabilizing solution such that the matrix
$A_{11}+B_{11}Z_2(i\infty)^{-1}(Q_{21}B_{11}+C_{11}^\dagger J_{11})^\dagger
$ is Hurwitz. The eigenvalues of this matrix are precisely the zeros of
$\tilde H_{21}(s)$~\cite{CALM-1997}, therefore $\tilde H_{21}^{-1}$ is
analytic in the half-plane $\mathrm{R}>-\tau$ ($\exists\tau>0)$, as
required in Step 5 of Algorithm~\ref{Alg1}.
 
\begin{remark}
As we remarked previously, the frequency domain condition~(\ref{eq:21})
which guarantees the contractiveness property~(\ref{eq:1}) is only a sufficient
condition. One can still attempt to apply the above procedure even
when~(\ref{eq:21}) fails. In this case, the properties 
$Z_1(i\infty)>0$, $Z_2(i\infty)>0$
can be ensured by expanding the set of the matrix inequalities to be solved
in Step 3 of Algorithm~\ref{Alg1} to include the following LMI in addition
to~(\ref{eq:40}):  
\begin{eqnarray}
\lefteqn{  I_{n+n_y}+
  \left[
  \begin{array}{cc}
    0_{n\times m_{11}} & I_n \\
    0_{n_y\times m_{11}} & 0_{n_y\times n}
  \end{array}
  \right]K^\dagger   \left[
  \begin{array}{cc}
    0_{m_{11}\times n} & 0_{m_{11}\times n_y} \\
    0_{n_y\times n} & I_n 
  \end{array}
  \right]} \nonumber && \\
&&
+ \left[\begin{array}{cc}
    0_{n\times m_{11}} & 0_{n\times n_y} \\
    0_{n_y\times m_{11}} & I_n 
  \end{array}
  \right]K
  \left[
  \begin{array}{cc}
    0_{m_{11}\times n} & 0_{m_{11}\times n_y} \\
    I_n & 0_{n\times n_y}
  \end{array}
  \right]>0. \nonumber \\
  \label{eq:43}
\end{eqnarray}
Indeed the left-hand side of~(\ref{eq:43}) is equal to $
\left[
  \begin{array}{cc}
    I_n & J_{11} \\
    J_{11}^\dagger & I_{n_y}
  \end{array}
\right]$, and so~(\ref{eq:43}) is equivalent to $Z_1(i\infty)>0$,
$Z_2(i\infty)>0$. 
\hfill$\Box$
\end{remark}

\section{Examples}\label{examples}

\subsection{Equalization of an optical cavity system}\label{example1}


 To illustrate the synthesis procedure described in the previous sections
 and demonstrate its practical applicability,
 we consider the quantum optical equalization system 
 shown in Fig.~\ref{cavity}. The channel 
  consists of an optical cavity and three optical beam splitters. It is
  similar to the optical cavity system considered in~\cite{UJ2}, however it
  includes an additional beam splitter. This additional element will allow
  us to illustrate application of our synthesis method in a low SNR
  scenario. Indeed, the effect of the environment noise $w_2$ of this beam
  splitter is similar to the effect of the additive Gaussian white noise in
  a conventional communication channel~\cite{HSK-1999}. As in the classical
  case, the low SNR situation occurs when $\sigma_{w_2}^2$ is sufficiently
  large.
\begin{figure}[t]
  \begin{center}
\psfragfig[width=\columnwidth]{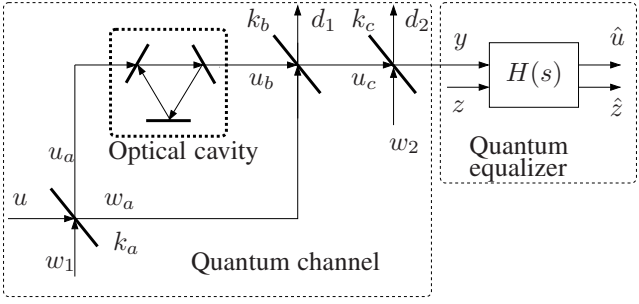}{ 
  \psfrag{Quantum channel}{Quantum channel}
  \psfrag{Quantum}{Quantum}
  \psfrag{equalizer}{equalizer}
  \psfrag{Optical cavity}{Optical cavity}
  \psfrag{b}{\hspace{-1ex}$w_1$}
  \psfrag{v}{$u$}
  \psfrag{bp}{$u_a$}
  \psfrag{vp}{$w_a$}
  \psfrag{bp2}{$u_b$}
  \psfrag{vp2}{}
  \psfrag{vp3}{$d_1$}
  \psfrag{alpha}{$k_a$}
  \psfrag{beta}{$k_b$}
  \psfrag{eta}{$k_c$}
  \psfrag{w}{$w_2$}
  \psfrag{z}{$z$}
  \psfrag{bhat}{${\hat u}$}
  \psfrag{wout}{$y_w$}
  \psfrag{vout}{${d}_2$}
  \psfrag{bout}{${u}_c$}
  \psfrag{zhat}{${\hat z}$}
  \psfrag{y}{${y}$}
  \psfrag{H}{\hspace{-2ex}$H(s)$}}
  \caption{A cavity, beam splitters and an equalizer system.}
  \label{cavity}
\end{center}
\end{figure}

The input
  operator $u$ in this example is scalar, and $\Sigma_u$ is a
  real constant. To emphasize this, we use the notation
  $\Sigma_u=\sigma_u^2$. The environment is represented by the quantum noises
  $w_1$, $w_2$, thus
  $w=\mathrm{col}(w_1,w_2)$ consists of two scalar operators. We assume that $\Sigma_w=
  \left[
    \begin{array}{cc}
     \sigma_{w_1}^2 & 0 \\
0 & \sigma_{w_2}^2 
    \end{array}
  \right]$. For simplicity, we assume that the
  transmittance parameters $k_a^2$, $k_b^2$, $k_c^2$ of the beam
  splitters are real positive numbers and that   
  $k_a$, $k_b$ and $k_c$ are also real positive constants, and that
  $k_a=k_b=k$. Thus, the relations between the inputs and
    outputs of the beam splitters are
    \begin{eqnarray*}
     && \left[
        \begin{array}{c}
   u_a \\ w_a       
        \end{array}
      \right]=
      \left[
        \begin{array}{cc}
          k& l \\ -l & k
        \end{array}
      \right]     
      \left[
        \begin{array}{c}
   u \\ w_1       
        \end{array}
      \right], \quad
\left[
        \begin{array}{c}
   u_c \\ d_1       
        \end{array}
      \right]=
      \left[
        \begin{array}{cc}
          k& l \\ -l & k
        \end{array}
      \right]     
      \left[
        \begin{array}{c}
   u_b \\ w_a       
        \end{array}
      \right], \\ 
&&   \left[
        \begin{array}{c}
   y \\ d_2       
        \end{array}
      \right]=
      \left[
        \begin{array}{cc}
          k_c& l_c \\ -l_c & k_c
        \end{array}
      \right]     
      \left[
        \begin{array}{c}
   u_c \\ w_2      
        \end{array}
      \right],  
    \end{eqnarray*}
where $l\triangleq \sqrt{1-k^2}$, $l_c\triangleq \sqrt{1-k_c^2}$ are real. 
    The  transfer function of the optical cavity is
$G_c(s)=\frac{s-\kappa+i\Omega}{s+\kappa+i\Omega}$, i.e., $u_b=G_c(s)u_a$;
$\kappa>0$, $\Omega$ are real numbers.  
Then the elements 
  of the transfer function $G(s)$ of the channel are  
  \begin{eqnarray}
    \label{eq:29}
G_{11}(s)&=&k_c(k^2G_c(s)-(1-k^2)), \nonumber \\
G_{12}(s)&=&
             \left[
             \begin{array}{cc}
k_ck\sqrt{1-k^2}(G_c(s)+1) &~~ \sqrt{1-k_c^2}
             \end{array}
              \right], \nonumber \\
G_{21}(s)&=&\left[
             \begin{array}{c}-k\sqrt{1-k^2}(G_c(s)+1)\\
-\sqrt{1-k_c^2}(k^2G_c(s)-(1-k^2))\end{array}
              \right], \nonumber \\
G_{22}(s)&=&\left[\begin{array}{cc}k^2-(1-k^2)G_c(s) & 0 \\
-\sqrt{1-k_c^2}k\sqrt{1-k^2}(G_c(s)+1) &~~ k_c\end{array}
              \right].
  \end{eqnarray}
Suppose $\sigma_{w_1}^2>\sigma_u^2>0$ and
$k^2<\frac{1}{2}$.  
Under these assumptions,
\begin{eqnarray}
&& \rho\triangleq
   1+\frac{\sigma_u^2}{2(\sigma_{w_1}^2-\sigma_u^2)k^2(1-k^2)}>1, \quad
   \hat\rho\triangleq \frac{\rho-1}{\rho+1}\in(0,1), \nonumber \\
&&  \delta\triangleq \frac{\sqrt{1-k^2}}{k}>1, \quad \hat\delta\triangleq
   \frac{\delta^2+1}{\delta^2-1}=\frac{1}{1-2k^2}>1.  
\label{eq:10} 
\end{eqnarray}
Using these notations, the function $\Psi(s)$ given in
equation~(\ref{eq:47}) is expressed as
\begin{eqnarray}
\label{eq:58}
  \Psi(s)&=&k_c^2\mu^2\frac{(s+i\Omega)^2-\hat\rho\kappa^2}{(s+i\Omega)^2-\kappa^2}+(1-k_c^2)\sigma_{w_2}^2,
\end{eqnarray}
where $\mu \triangleq \sqrt{2(\sigma_{w_1}^2-\sigma_u^2)k^2(1-k^2)(1+\rho)}$.
This gives the expression for the matrix $\Phi(s)$ in
equation~(\ref{eq:75}),
\begin{equation}
  \label{eq:67}
  \Phi(s)=
  \left[
    \begin{array}{cc}
      k_c^2\mu^2\frac{(s+i\Omega)^2-\hat\rho\kappa^2}{(s+i\Omega)^2-\kappa^2}+(1-k_c^2)\sigma_{w_2}^2~&~
k_c\frac{\sigma_u^2+1}{\hat\delta}\frac{s+\hat\delta\kappa+i\Omega}{s+\kappa+i\Omega}\\
k_c\frac{\sigma_u^2+1}{\hat\delta}\frac{s-\hat\delta\kappa+i\Omega}{s-\kappa+i\Omega}
  & \sigma_u^2+2 
    \end{array}
  \right].
\end{equation}

\begin{proposition}\label{prop.rho}
If
\begin{equation}
  \label{eq:71}
  \lambda^2>\frac{k_c^2(1+\sigma_u^2)^2}{k_c^2\mu^2\hat\rho+(1-k_c^2)\sigma_{w_2}^2}-(\sigma_u^2+2),
\end{equation}
then Assumption~\ref{A.rho} is satisfied for the considered system. 
\end{proposition}

\emph{Proof: } The proof of the proposition is by direct construction. 
Condition~(\ref{eq:71}) guarantees that
\begin{eqnarray*}
  \mu_1&=&(k_c^2\mu^2\hat\rho+(1-k_c^2)\sigma_{w_2}^2)(\sigma_u^2+2+\lambda^2) \\
&& -k_c^2(1+\sigma_u^2)^2>0, \\
  \mu_2&=&(k_c^2\mu^2+(1-k_c^2)\sigma_{w_2}^2)\hat\delta^2(\sigma_u^2+2+\lambda^2)\\
&&-k_c^2(1+\sigma_u^2)^2>0.
\end{eqnarray*}
The first inequality follows directly from~(\ref{eq:71}). To obtain the second
inequality, recall that $\hat\delta^2>1>\hat\rho$. Therefore 
\begin{eqnarray*}
  (k_c^2\mu^2+(1-k_c^2)\sigma_{w_2}^2)\hat\delta^2> 
k_c^2\mu^2\hat\rho+(1-k_c^2)\sigma_{w_2}^2,
\end{eqnarray*}
hence $\mu_2>\mu_1>0$. 
Therefore, the following constants are well
defined, 
\begin{eqnarray*}
  \label{eq:80a1}
  \alpha_1=\left(
\frac{\mu_2}{\hat\delta^2(\sigma_u^2+2+\lambda^2)}
\right)^{1/2}, \quad 
  \beta_1=\left(
\frac{\mu_1}{\mu_2}
\right)^{1/2}.
\end{eqnarray*}
Also, let
\begin{eqnarray*}
  \label{eq:80a2}
  \alpha_2&=&
\frac{k_c(1+\sigma_u^2)}{\hat\delta(\sigma_u^2+2+\lambda^2)^{1/2}}.
\end{eqnarray*}
It is readily
verified that the transfer function 
\begin{eqnarray}
  \label{eq:77}
\Upsilon_\lambda(s)
=
  \left[
  \begin{array}{cc}
\alpha_1\frac{s+\beta_1\hat\delta\kappa+i\Omega}{s+\kappa+i\Omega}
    &\alpha_2\frac{s+\hat\delta\kappa+i\Omega}{s+\kappa+i\Omega} \\
0 & \sqrt{\sigma_u^2+2+\lambda^2}
  \end{array}
  \right]
\end{eqnarray}
is a spectral factor of $\Phi_\lambda$\footnote{The choice of
the spectral factor $\Upsilon_\lambda(s)$ is unique up to a unitary
transformation; see~\cite{Anderson-1967,Youla-1961} for details.}.  
\hfill$\Box$

\begin{remark}
Of course, one can use the results in~\cite{Anderson-1967} to check whether
$\Phi_\lambda(s)$ in this example admits the required spectral
factorization. This will involve checking whether a certain algebraic Riccati
equation has a Hermitian solution. In this example, this Riccati
equation reduces to a scalar quadratic equation, and
condition~(\ref{eq:71}) is necessary and 
sufficient for this quadratic equation to have a real scalar
solution. This can be checked directly, although the calculations are
unwieldy to include here.    
\hfill$\Box$
\end{remark}

The transfer function $\Upsilon_\lambda(s)$ in~(\ref{eq:77}) can be 
expressed as
\begin{eqnarray}
  \label{eq:83}
\Upsilon_\lambda(s)&=&
\left[
  \begin{array}{c}
    1 \\ 0
  \end{array}
\right](s+\kappa+i\Omega)^{-1}\left[
  \begin{array}{cc}
    \alpha_1(\beta_1\hat\delta-1)\kappa~ &~
                                      \alpha_2(\hat\delta-1)\kappa 
  \end{array}
\right] \nonumber \\
&+&
\left[
  \begin{array}{cc}
\alpha_1 & \alpha_2 \\
0 & \sqrt{\sigma_u^2+2+\lambda^2} 
  \end{array}
\right].
\end{eqnarray}
According to this, we choose the following coefficients for the
state-space model of $\Upsilon_\lambda(s)$: 
\begin{eqnarray}
  \label{eq:84}
  A_\lambda&=&-(\kappa+i\Omega), \quad B_\lambda=\left[
     \begin{array}{cc}
 \alpha_1(\beta_1\hat\delta-1)\kappa~ &~ \alpha_2(\hat\delta-1)\kappa
     \end{array}
   \right], \nonumber \\ 
C_{1,\lambda}&=&1, \quad C_{2,\lambda}=0, \nonumber \\
D_{1,\lambda}&=& \left[
     \begin{array}{cc}
\alpha_1~ &~ \alpha_2 
     \end{array}
   \right], \quad
D_{2,\lambda}= \left[
     \begin{array}{cc}
0~ &~ \sqrt{\sigma_u^2+2+\lambda^2} 
     \end{array}
   \right].
\end{eqnarray}

We now demonstrate the 
application of Algorithm~\ref{Alg1} to computing a transfer
function $H(s)$ of a coherent guaranteed cost equalizer. 
For this, we select the same numerical values for the parameters
$\sigma_u^2$, $\sigma_{w_1}^2$, $k$, $\kappa$, and $\Omega$ as those used 
in a similar example in~\cite{UJ2}: $\sigma_u^2=0.1$, $\sigma_{w_1}^2=0.2$,
$k=0.4$, $\kappa=5\times 10^8$, $\Omega=10^9$. Also, we let $k_c^2=0.5$, and selected
$\sigma_{w_2}^2=3$. 
 
In this example, condition~(\ref{eq:21}) requires that there must exist
a constant $\theta>0$ such that 
\begin{eqnarray}
  \label{eq:5}
\lefteqn{
  \frac{1}{\theta^2}+((\sigma_u^2+2-\gamma^2)-\Psi(i\omega))\frac{1}{\theta}} &&
  \nonumber \\
&&-\det\left(\Phi(i\omega)-
  \left[
  \begin{array}{cc}
    0 & 0 \\ 0 & \gamma^2
  \end{array}
  \right]\right)<0 \quad \forall \omega\in\mathbf{R}.
\end{eqnarray}
By decreasing $\gamma^2$ in small
steps starting from $\gamma^2=2+\sigma_u^2=2.1$, it was found that this
inequality holds with $\theta=1.2956$ for all
$\gamma^2\in(0,\gamma_0^2]$, where $\gamma_0^2=2.0937$. Also, with the chosen
numerical parameter values,  
the expression on the right hand side of inequality~(\ref{eq:71}) is
negative ($=-1.7097$), therefore we let $\lambda^2=0$.   

With these parameters, the optimization problem~(\ref{eq:38conv}) was solved
numerically subject to the additional constraint $\bar\gamma^2\le
\gamma_0^2+\lambda^2=2.0937$. According to Corollary~\ref{main.T.conv}, 
for this range of $\gamma^2$ the corresponding
transfer function~(\ref{eq:51}) is guaranteed to be contractive whenever the
constraints of the problem~(\ref{eq:38conv}) are feasible. The
optimization resulted in the optimal value 
$\bar\gamma^2_*\approx 1.9255$; this yields $\gamma_*^2\approx 1.9255\in
(0,2.0937]$ since $\lambda^2=0$. Next, we chose an acceptable level of
guaranteed mean-square equalization performance to be within a 1\%
margin of $\gamma_*^2$, i.e., we let    
$\gamma^2=1.01\gamma_*^2=1.9448\in(\gamma_*^2,2.0937]$. Such
level of accuracy is acceptable for the purpose of this illustration. 
According to
Theorem~\ref{SDP.primal.LMI} this value of $\gamma^2$ approximates the
theoretically optimal value $\gamma_\circ^2$ achievable by means of
completely passive filters within an error margin of 1\%. Then, we carried 
out Steps~1--3 of Algorithm~\ref{Alg1} with
$\bar\gamma^2=\gamma^2+\lambda^2=1.9448$. This produced the transfer function
\begin{eqnarray}
  \label{eq:56}
&&  H_{11}(s)=a\frac{s+b\kappa_1+i\Omega}{s+\kappa_1+i\Omega}, \\
&&  a=-0.36292, \quad b=0.9837, \quad \kappa_1=3.1853\times 10^8. \nonumber
\end{eqnarray}
The magnitude Bode plot of this transfer function shown in
Figure~\ref{fig:Bode} confirms that it is contractive.
\begin{figure}[t]
  \centering
  \includegraphics[width=0.85\columnwidth]{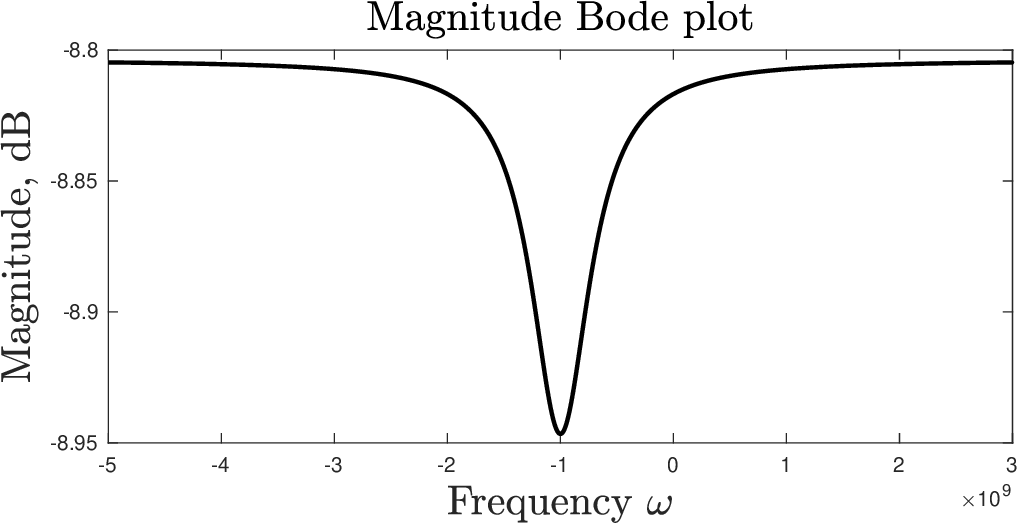}
  \caption{The magnitude Bode plot of the
suboptimal transfer function $H_{11}(s)$ in equation~(\ref{eq:56}).} 
  \label{fig:Bode}
\end{figure}
In addition, this $H_{11}(s)$
is stable and analytic in the half-plane $\mathrm{Re}s>-\kappa_1$. 
Thus, the corresponding set $\mathcal{H}_{11,\gamma}^-$ is not empty; this
confirms the result of Theorem~\ref{main.T}. According to
Theorem~\ref{P1a>P1}, a physically realizable
equalizer can be constructed from this 
$H_{11}(s)$. We constructed such equalizer following Steps 4 and 5 of
Algorithm~\ref{Alg1}. 

First, $Z_1(s)$ and $Z_2(s)$ in equation~(\ref{eq:26}) were computed. For
this example, we obtained   
\begin{eqnarray}
  \label{eq:28}
  Z_1(s)&=&Z_2(s) 
=(1-a^2)\frac{(s+i\Omega)^2-\frac{1-b^2a^2}{1-a^2}\kappa_1^2}{(s+i\Omega)^2-\kappa_1^2}. 
\end{eqnarray}
It is easy to check that $Z_1(s)$ and $Z_2(s)$ are para-Hermitian and 
$Z_1(i\omega)>0$, $Z_2(i\omega)>0$ in the closed imaginary axis if and only
if $a^2<1$, $a^2b^2<1$. Clearly, the parameter values in~(\ref{eq:56}) meet this
condition. The Riccati 
equations~(\ref{eq:23}), (\ref{eq:24}) take the form
\begin{eqnarray}
  \label{eq:49}
&&  -2\kappa_1
  Q_{12}+\frac{(Q_{12}+a^2(b-1)\kappa_1)^2}{1-a^2} \nonumber \\
&& \hspace{4.3cm} +a^2(b-1)^2\kappa_1^2=0, 
  \\
&&  -2\kappa_1 Q_{21}+\frac{a^2}{1-a^2}(Q_{21}(b-1)\kappa_1+1)^2+1=0.
\label{eq:50}
\end{eqnarray}
Under the conditions $a^2<1$, $a^2b^2<1$, their  
stabilizing positive solutions are, respectively,
\begin{eqnarray}
  \label{eq:54}
&&  Q_{12}=(1-a^2b-(1-a^2b^2)^{1/2}(1-a^2)^{1/2})\kappa_1, \\
&&  Q_{21}=\frac{1-a^2b-(1-a^2b^2)^{1/2}(1-a^2)^{1/2}}{(b-1)^2a^2\kappa_1}.
\end{eqnarray}
Then spectral
factors of $Z_1(s)$, $Z_2(s)$ were obtained
using~(\ref{eq:27}),~(\ref{eq:34}):  
\begin{eqnarray}
  \label{eq:30}
   H_{12}(s)&=&-\sqrt{1-a^2}\frac{s+\sqrt{\frac{1-b^2a^2}{1-a^2}}\kappa_1+i\Omega}{s+\kappa_1+i\Omega}, \nonumber
  \\
\tilde H_{21}(s)&=&-H_{12}(s)\nonumber \\
&=&\sqrt{1-a^2}\frac{s+\sqrt{\frac{1-b^2a^2}{1-a^2}}\kappa_1+i\Omega}{s+\kappa_1+i\Omega}.
\end{eqnarray}
Note that $\tilde
  H_{21}^{-1}(s)$ is stable and analytic in the right half-plane
  $\mathrm{Re}s>-\sqrt{\frac{1-b^2a^2}{1-a^2}}\kappa_1$ of
  the complex plane, as expected. Finally, to carry out Step 5 we selected 
\begin{equation*}
  \label{eq:63}
  U(s)=\frac{s-\sqrt{\frac{1-b^2a^2}{1-a^2}}\kappa_1+i\Omega}{s+\sqrt{\frac{1-b^2a^2}{1-a^2}}\kappa_1+i\Omega}. 
\end{equation*}
This transfer function is paraunitary, stable and analytic in the right
half-plane $\mathrm{Re}s>-\sqrt{\frac{1-b^2a^2}{1-a^2}}\kappa_1$ and also
satisfies all other requirements of the proposed procedure.  
Using it,  the remaining blocks of $H(s)$ were obtained
according to~(\ref{eq:81}):
\begin{eqnarray}
  \label{eq:72}
  H_{21}(s)&=&U(s)\tilde H_{21}(s) \nonumber \\
&=&\sqrt{1-a^2}\frac{s-\sqrt{\frac{1-b^2a^2}{1-a^2}}\kappa_1+i\Omega}{s+\kappa_1+i\Omega},
    \nonumber \\
  H_{22}(s)&=&-U(s)(\tilde H_{21}^{-1}(s))^HH_{11}(s)^HH_{12}(s) \nonumber \\
&=&a\frac{s-b\kappa_1+i\Omega}{s+\kappa_1+i\Omega}.
\end{eqnarray}

The plot of the error power spectrum density $P_e(i\omega)$ for the
resulting equalizer shown in Fig.~\ref{fig.compare} confirms that with this
equalizer,  
$\sup_\omega P_e(i\omega)<\gamma^2$. 
\begin{figure}[t]
  \centering
  \includegraphics[width=0.95\columnwidth]{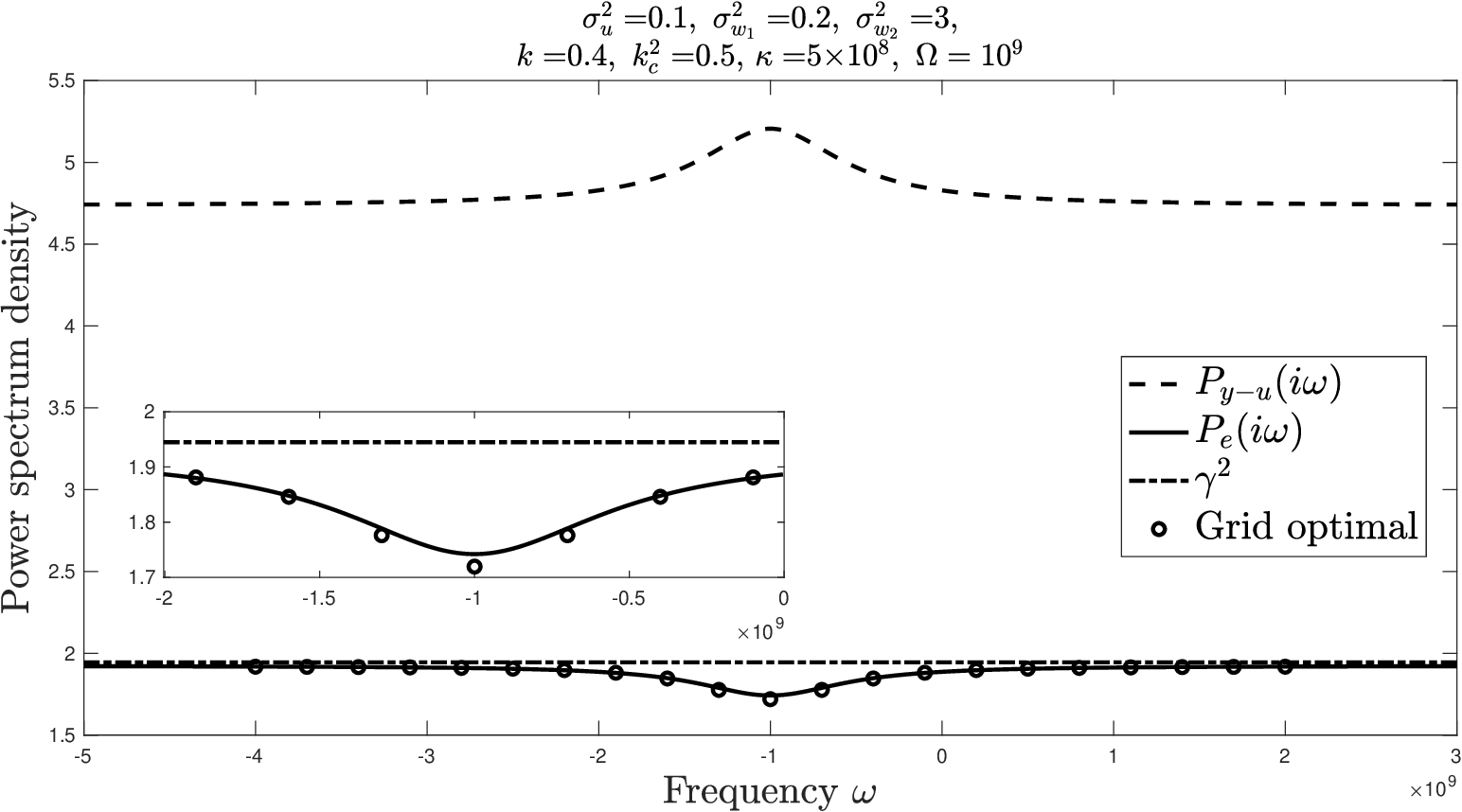}
  \caption{Power spectrum densities $P_{y-u}$ and $P_{e}$ and the suboptimal
    value $\gamma^2$. The circles indicate the values 
      $P_e(i\omega_l)$ obtained 
    from the solutions of the optimization problem~(\ref{eq:71.LMI}).} 
  \label{fig.compare}
\end{figure}
Also, Fig.~\ref{fig.compare} compares the power spectrum density $P_e$ with
the power spectrum density of the difference between the channel input and
output $y-u$. The figure confirms that the output of the equalizer
represents $u$ with a substantially higher fidelity (in the mean-square
sense) then the channel output $y$. 

Also, Fig.~\ref{fig.compare} compares the proposed method with
  the method developed 
  in~\cite{UJ2}. The method allows to construct a physically realizable
equalizer which minimizes the power spectrum density $P_e(i\omega)$ at
selected frequency points $\omega_l$. $l=1,\ldots,L$. The equalizer is
obtained by first solving the semidefinite program 
\begin{eqnarray}
  \label{eq:71.LMI}
  && \nu_\circ^2=\inf \nu^2 \\
  && \mbox{s.t. }     \left[
     \begin{array}{cc}
\Xi_{11,l} & H_{11,l}M(i\omega_l) \\
M(i\omega_l)^\dagger H_{11,l}^\dagger & -I_q 
     \end{array}
     \right]< 0, \nonumber \quad \\
  && \phantom{\mbox{s.t. }} 
     \left[
     \begin{array}{cc}
       I_n & H_{11,l}\\ H_{11,l}^\dagger & I_n
     \end{array}
     \right]\ge 0 \quad \forall l=1,\ldots,L, \nonumber
\end{eqnarray}
in variables $\nu^2$ and $H_{11,l}$. Here,
\begin{eqnarray*}
\Xi_{11,l}&\triangleq&
(2-\nu^2)I+\Sigma_u^T-H_{11,l}G_{11}(i\omega_l)(I+\Sigma_u^T) \\
&& -(I+\Sigma_u^T)G_{11}(i\omega_l)^\dagger H_{11,l}^\dagger,  
\end{eqnarray*}
$M(s)$ is a stable $n_y\times q$ spectral factor of $\Psi(s)$ which is
analytic in the complex half-plane $\mathrm{Re}\,s>-\tau$ ($\exists \tau>0$):
\[
\Psi(s)=M(s)M(s)^H.
\]
The obtained $\{H_{11,l},~l=1,\ldots,L\}$, are then used as the data for
Nevanlinna-Pick interpolation. The interpolation produces a transfer
function $H_{11}(s)$ such that $H_{11}(i\omega)H_{11}(i\omega)^\dagger\le
I_n$ while $H_{11}(i\omega_l)=H_{11,l}$. This transfer
function is shown in~\cite{UJ2} to have the stability, causality and analiticity
properties which allow one to construct a paraunitary transfer function
$H(s)$ for which the inequality (\ref{eq:6'.sub}) is satisfied at the
selected points $\omega_l$ with $\gamma^2=\nu^2$, i.e.,  
\begin{equation}
  P_e(\omega_l)< \nu^2 I, \quad l=1,\ldots,L,
\label{eq:41}
\end{equation}
for any feasible value $\nu^2$ of the problem~(\ref{eq:71.LMI}).
From~(\ref{eq:71.LMI}) and (\ref{eq:6'.sub}), it is clear that
$\gamma_\circ^2\ge \sup \nu_\circ^2$ where the supremum is over all
possible grids $\{\omega_l, l=1,\ldots,L\}$. 

We applied this method to compute the optimal $\nu_\circ^2=1.9191$ and the
corresponding values $P_e(i\omega_l)$ at 21 frequency points $\omega_l$ evenly
distributed in the interval $[-4\times 10^9,2\times 10^9]$. These
values are shown as circles in Fig.~\ref{fig.compare}. From this
comparison we conclude that in our example, the equalization method
developed in this paper produces a mean-square error whose PSD
$P_e(i\omega)$ is indeed near optimal. However, the 
interpolation technique from~\cite{UJ2} does not guarantee that the PSD
$P_e(i\omega)$ is less than $\gamma^2I$ at frequencies other than the
selected grid points $\omega_l$, $l=1,\ldots,L$. In contrast, the method in this
paper guarantees the near optimal performance over the entire range of
frequencies. 
 
We conclude this example by pointing out that the equalizer found in
this example can be realized as an interconnection of an optical cavity and two
beam splitters, as shown in Figure~\ref{fig:cavity.filter}. 
The optical cavity used in this realization is
\begin{equation*}
y_2=H_c(s)y_1, \quad  H_c(s)=\frac{s-\kappa_1+i\Omega}{s+\kappa_1+i\Omega}.
\end{equation*}
The beam splitters' operators are
    \begin{eqnarray*}
     && \left[
        \begin{array}{c}
   y_1 \\ z_1       
        \end{array}
      \right]=
      \left[
        \begin{array}{cc}
          \xi_1 & \eta_1 \\ \eta_1 & -\xi_1
        \end{array}
      \right]     
      \left[
        \begin{array}{c}
   y \\ z       
        \end{array}
      \right], \quad
     \left[
        \begin{array}{c}
   \hat u \\ \hat z       
        \end{array}
      \right]=
      \left[
        \begin{array}{cc}
          \eta_2& \xi_2 \\ \xi_2 & -\eta_2
        \end{array}
      \right]     
      \left[
        \begin{array}{c}
   y_2 \\ z_2       
        \end{array}
      \right],  
    \end{eqnarray*}
where
\begin{eqnarray*}
 && \eta_1=-\sqrt{\frac{1+a^2b-\sqrt{(1-a^2b^2)(1-a^2)}}{2}}, \\
 && \xi_1=\sqrt{1-\eta_1^2}=\sqrt{\frac{1-a^2b+\sqrt{(1-a^2b^2)(1-a^2)}}{2}}, \\
 && \eta_2=-\sqrt{\frac{1-a^2b-\sqrt{(1-a^2b^2)(1-a^2)}}{2}}, \\
 && \xi_2=\sqrt{1-\eta_2^2}=\sqrt{\frac{1+a^2b+\sqrt{(1-a^2b^2)(1-a^2)}}{2}},
\end{eqnarray*}
and $a$, $b$ are the numerical constants defined in~(\ref{eq:56}). 
\begin{figure}[t]
  \begin{center}
\psfragfig[width=0.7\columnwidth]{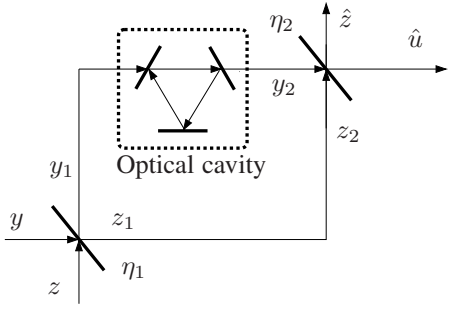}{ 
  \psfrag{Optical cavity}{Optical cavity}
  \psfrag{b}{\hspace{-1ex}$z$}
  \psfrag{v}{$y$}
  \psfrag{bp}{$y_1$}
  \psfrag{vp}{$z_1$}
  \psfrag{bp2}{$y_2$}
  \psfrag{vp2}{$z_2$}
  \psfrag{alpha}{$\eta_1$}
  \psfrag{beta}{$\eta_2$}
  \psfrag{w}{\hspace{-3ex}$w$}
  \psfrag{bhat}{$\hat u$}
  \psfrag{wout}{$y_w$}
  \psfrag{vout}{$\hat z$}
  \psfrag{y}{${\hat u}$}
  \psfrag{H}{\hspace{-2ex}$H(s)$}}
  \caption{A cavity and beam splitters realization of the equalizer.}
  \label{fig:cavity.filter}
\end{center}
\end{figure}

It is worth noting that in this example, when $\sigma_{w_2}^2$ is reduced,
condition~(\ref{eq:21}) will fail eventually for all $\gamma>\gamma_*$. Despite
this, the transfer function $H_{11}(s)$ obtained using Algorithm~\ref{Alg1}
remains contractive for 
$\sigma_{w_2}^2\ge 0.975$; see Fig.~\ref{fig:Bode-1}. This
confirms that Theorem~\ref{SDP.primal.LMI} is sufficient but it is not
necessary to guarantee the contractiveness property~(\ref{eq:1}). 
      
\begin{figure}[t]
  \centering
  \includegraphics[width=0.85\columnwidth]{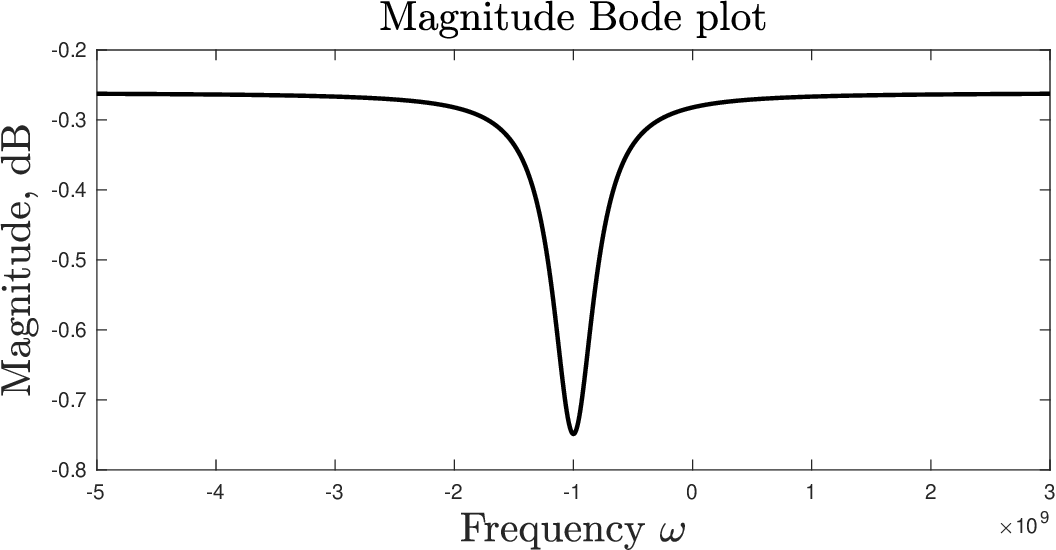}
  \caption{The magnitude Bode plot of the
suboptimal transfer function $H_{11}(s)$ for
$\sigma_{w_2}^2=0.975$. Inequality~(\ref{eq:21}) fails to hold for all
$\gamma\in[\gamma_*,2+\sigma_u^2]$. Despite this, $\|H_{11}\|_\infty<1$.} 
  \label{fig:Bode-1}
\end{figure}

\subsection{Equalization of a two-input two-output four-noise system
  comprised of two cavities}\label{example2} 

%
\begin{figure}[t]
  \begin{center}
\psfragfig[width=0.9\columnwidth]{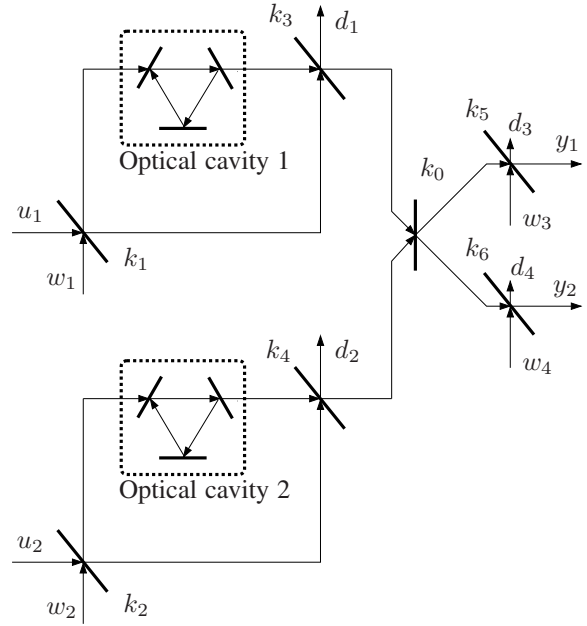}{ 
  \psfrag{Optical cavity 1}{Optical cavity 1}
  \psfrag{Optical cavity 2}{Optical cavity 2}
  \psfrag{w1}{\hspace{-1.5ex}$w_1$}
  \psfrag{w2}{\hspace{-1.5ex}$w_2$}
  \psfrag{w3}{$w_3$}
  \psfrag{w4}{$w_4$}
  \psfrag{u1}{$u_1$}
  \psfrag{u2}{$u_2$}
  \psfrag{d1}{$d_1$}
  \psfrag{d2}{$d_2$}
  \psfrag{d3}{$d_3$}
  \psfrag{d4}{$d_4$}
  \psfrag{k1}{$k_1$}
  \psfrag{k2}{$k_2$}
  \psfrag{k3}{$k_3$}
  \psfrag{k4}{$k_4$}
  \psfrag{k0}{$k_0$}
  \psfrag{k5}{$k_5$}
  \psfrag{k6}{$k_6$}
  \psfrag{y1}{$y_1$}
  \psfrag{y2}{$y_2$}
}
   \caption{A two-input two-output four-noise system comprised of two
     cavities and seven beam splitters.}
   \label{two-cavity-system}
\end{center}
\end{figure}
We now illustrate the proposed synthesis method using a system comprised of
an interconnection of two quantum optical systems considered in the
previous section. The system is shown in Fig.~\ref{two-cavity-system}. The
system has two input fields indicated as $u_1$, $u_2$ in the figure. Also,
it interacts with four environment fields, designated in the figure as the
inputs $w_1$, $w_2$, $w_3$ and $w_4$. The transfer functions of the optical
cavities are
\[
G_{c,l}(s)=\frac{s-\kappa_l+i\Omega_l}{s+\kappa_l+i\Omega_l}, \quad l=1,2.
\] 
The output fields to be equalized are $y_1$, $y_2$, while the outputs $d_1$,
$d_2$, $d_3$, $d_4$ describe losses to the environment. The task in this
section is to obtain a physically realizable transfer function $H(s)$ which
takes the optical field $y=\mathrm{col}(y_1,y_2)$ as its
input (along with some vacuum noise) and produces `optimally 
equalized' optical field $\hat u=\mathrm{col}(\hat u_1,\hat u_2)$ as its
output.    

From Fig.~\ref{two-cavity-system}, the components $G_{11}(s)$, $G_{12}(s)$
of the system transfer function $G(s)$ are
\begin{eqnarray*}
  \label{eq:64}
&& G_{11}(s)=
  \left[
  \begin{array}{cc}
   a_{11}\frac{s-\frac{\kappa_1}{2k_1^2-1}+i\Omega_1}{s+\kappa_1+i\Omega_1} &
    a_{12}\frac{s-\frac{\kappa_2}{2k_2^2-1}+i\Omega_2}{s+\kappa_2+i\Omega_2} \\
    a_{21}\frac{s-\frac{\kappa_1}{2k_1^2-1}+i\Omega_1}{s+\kappa_1+i\Omega_1} &
    a_{22}\frac{s-\frac{\kappa_2}{2k_2^2-1}+i\Omega_2}{s+\kappa_2+i\Omega_2} \\
  \end{array}
  \right],\\
&& G_{12}(s)=
  \left[
  \begin{array}{cccc}
   b_{11}\frac{s+i\Omega_1}{s+\kappa_1+i\Omega_1} &
    b_{12}\frac{s+i\Omega_2}{s+\kappa_2+i\Omega_2} & \sqrt{1-k_5^2} & 0 \\
    b_{21}\frac{s+i\Omega_1}{s+\kappa_1+i\Omega_1} &
    b_{22}\frac{s+i\Omega_2}{s+\kappa_2+i\Omega_2} & 0 & \sqrt{1-k_6^2} 
  \end{array}
  \right].\label{eq:65}
\end{eqnarray*}
Here we used the notation
\begin{eqnarray*}
&& a_{11}= k_5k_0(2k_1^2-1), \quad a_{12}=k_5\sqrt{1-k_0^2}(2k_2^2-1),
   \nonumber \\
&& a_{21}=-k_6\sqrt{1-k_0^2}(2k_1^2-1), \quad a_{22}=k_6k_0(2k_2^2-1),
   \nonumber \\
&& b_{11}=2k_5k_0k_1\sqrt{1-k_1^2}, \quad
   b_{12}=2k_5\sqrt{1-k_0^2}k_2\sqrt{1-k_2^2}, \nonumber \\
&& b_{21}=-2k_6\sqrt{1-k_0^2}k_1\sqrt{1-k_1^2}, \quad
   b_{22}=2k_6k_0k_2\sqrt{1-k_2^2}.
\end{eqnarray*}
Using these expressions, closed form expressions for
$\Phi_\lambda(s)$ and $\Upsilon_\lambda(s)$ can be derived. However, the
derivation is quite tedious. Therefore in this example we carried out the
synthesis of an equalizer numerically.     

As in the previous example, we use equal numerical values for the
parameters $k_1$, $k_3$
and $k_2$, $k_4$ of the beam splitters: $k_1=k_3=0.4$, $k_2=k_4=0.3$. Also,
we let 
$k_0=k_5=k_6=1/\sqrt{2}$. The parameters of the cavities were chosen as
follows: $\kappa_1=7.5\times 10^8$, $\kappa_2=3\times 10^8$,
$\Omega_1=10^9$, $\Omega_2=-0.5\times 10^9$.

The input fields $u_1$, $u_2$ $w_1$, and $W_2$ were assumed to be
uncorrelated, however the fields $w_3$, $w_4$ were allowed
to be correlated, 
\[
\Sigma_u=
\left[
  \begin{array}{cc}
    \sigma_{u_1}^2& 0 \\ 0 & \sigma_{u_2}^2
  \end{array}
\right], \quad
\Sigma_w=
\left[
  \begin{array}{cccc}
    \sigma_{w_1}^2& 0 & 0 & 0\\ 
    0 & \sigma_{w_2}^2 & 0 & 0 \\
    0 & 0 & \sigma_{w_3}^2 & \sigma_{w_3,w_4} \\
    0 & 0 & \sigma_{w_3,w_4} & \sigma_{w_4}^2
  \end{array}
\right].
\]    
The numerical values for the matrices $\Sigma_u$, $\Sigma_w$ were selected
as follows: 
$\sigma_{u_1}^2=0.1$, $\sigma_{u_2}^2=0.2$, $\sigma_{w_1}^2=0.2$,
$\sigma_{w_2}^2=0.3$, $\sigma_{w_3}^2=\sigma_{w_4}^2=3$, $\sigma_{w_3,w_4}=0.2$.   
With these parameters and $\lambda^2=0$, the minimum eigenvalue of the
matrix $\Phi_\lambda(i\omega)$ was computed across a range of frequencies
and was found to be greater than $0.9320$. That is,
Assumption~\ref{A.rho} is satisfied with $\lambda^2=0$. This allowed us to
compute the spectral factor $\Upsilon_\lambda(s)$. For this, a spectral
factor of $\Phi_\lambda(s)$ was found using the Matlab function
\texttt{spectralfact}, then its unobservable and 
uncontrollable states were removed. 

Then the optimization problem in Step 1 of
Algorithm~\ref{Alg1} was solved. As in the previous example, we
used the convex problem~(\ref{eq:38conv}) at this step, and the
dimension of the state space model for the filter $H_{11}$ was set to be
equal to the dimension of the state space model of $\Upsilon_\lambda$,
$m_{11}=m=2$.  
The value $\bar\gamma_*^2$ of the optimization problem~(\ref{eq:38conv})
$\bar\gamma_*^2$ was computed to be $\bar\gamma_*^2\approx 1.9209$. 

Next, a near optimal value $\gamma^2= 1.9401>\bar\gamma_*^2-\lambda^2\approx
1.9209$ was selected. It was found that condition~(\ref{eq:21}) of
Theorem~\ref{SDP.primal.LMI} is satisfied with this $\gamma^2$ and
$\theta=1.5$ (the eigenvalues of the matrix on the left hand side
of condition~(\ref{eq:21}) were found to be greater than 0.0268 across a
range of frequencies).    
According to Theorem~\ref{SDP.primal.LMI}, it follows that 
$\gamma_\circ^2=(\gamma_\circ')^2=\gamma_\circ''=\gamma_*^2$. That is, the
tightest bound on the equalization performance achievable by the proposed
method is $\gamma_\circ^2\approx 1.9209$. Also, the set
$\mathcal{H}_{11,\gamma}^-$ corresponding to the chosen near optimal
$\gamma^2=1.9401$ is not empty, according to 
Theorem~\ref{main.T}. Following Steps 2 and 3 of Algorithm~\ref{Alg1} an
element $H_{11}(s)$ of this set was computed. The parameters of the
obtained $2\times 2$ transfer function matrix $H_{11}(s)$ are presented in
the Appendix. It was verified that $\|H_{11}\|_\infty=0.4054$, i.e., the obtained
$H_{11}(s)$ is indeed contractive.

Next, Steps 4 and 5 of Algorithm~\ref{Alg1} were carried out to
obtain the remaining components of the equalizing filter. 
First, we confirmed that the matrices $Z_1(i\omega)$ and $Z_2(i\omega)$
were full rank matrices. This was done by computing the minimum eigenvalue
of these matrices across a range of frequencies, which was found to be
greater than $0.8357$.
Thus, we let $r=n=2$, and let the equalizer  $H(s)$ in this example be a
$4\times 4$ transfer function comprised of the $2\times 2$ transfer function
$H_{11}(s)$ found above, and $2\times 2$ transfer functions $H_{12}$,
$H_{21}(s)$ and $H_{22}(s)$. To compute these transfer functions, the
stabilizing  solutions of the Riccati
equations~(\ref{eq:23}),~(\ref{eq:24}) were obtained. Then the transfer
functions $H_{12}(s)$, $\tilde H_{21}(s)$ were obtained using 
equations~(\ref{eq:27}),~(\ref{eq:34}). The 
resulting transfer function $\tilde H_{21}^{-1}(s)$ was found to be stable and
analytic in the closed right half-plane of the complex plane; 
its poles were found to be located at $(-0.7126 - 1.0000i)\times 10^9$ and $
(-0.3236+0.5000i) \times 10^9$. Accordingly, all poles of
the transfer function $(\tilde H_{21}^{-1}(s))^HH_{11}(s)^H$ were found to
lie in the right half-plane of the complex plane.
To cancel these poles when constructing $H_{21}$ and $H_{22}$, first a
minimal realization of $(\tilde H_{21}^{-1}(s))^HH_{11}(s)^H$ was computed.
Its state matrix was found to have the eigenvalues 
\begin{eqnarray*}
  && p_1=(0.7126 -  1.0000i)\times 10^9, \\
  && p_2=(0.3236+0.5000i)\times 10^9,
\end{eqnarray*}
These values were used to construct the 2nd-order scalar unitary transfer
function $U(s)$, 
\[
U(s)=\frac{(s-p_1)(s-p_2)}{(s+p_1^*)(s+p_2^*)},
\]
to be used in equations~(\ref{eq:81}). 
 
With this choice of $U(s)$, the remaining blocks of the equalizer $H(s)$
were computed from equations~(\ref{eq:81}). Their values are presented in
the Appendix. With the computed blocks, we formed the matrix $H(s)$ and
found that $\|H(s)H(s)^H-I\|_\infty=7.64132\times
10^{-14}$. This shows that the computed equalizer transfer function $H(s)$ is
paraunitary with a high numerical precision.  

The graphs of the largest eigenvalues of the PSD matrices $P_e(i\omega)$
corresponding to the found equalizer $H(s)$ and $P_{y-u}(i\omega)$ are
shown in Fig.~\ref{fig:two-cavities-PSDs}. The graphs indicate that the
obtained equalizer $H(s)$ reduces the mean-square error 
between $u$ and $\hat u$ to a value below the set near optimal
$\gamma^2$. This mean square error is lower than the mean-square errors
between $y$ and $u$ and between $-y$ and $u$; the latter is the
error relative the phase inverted channel output.       
\begin{figure}[t]
  \centering
  \includegraphics[width=0.85\columnwidth]{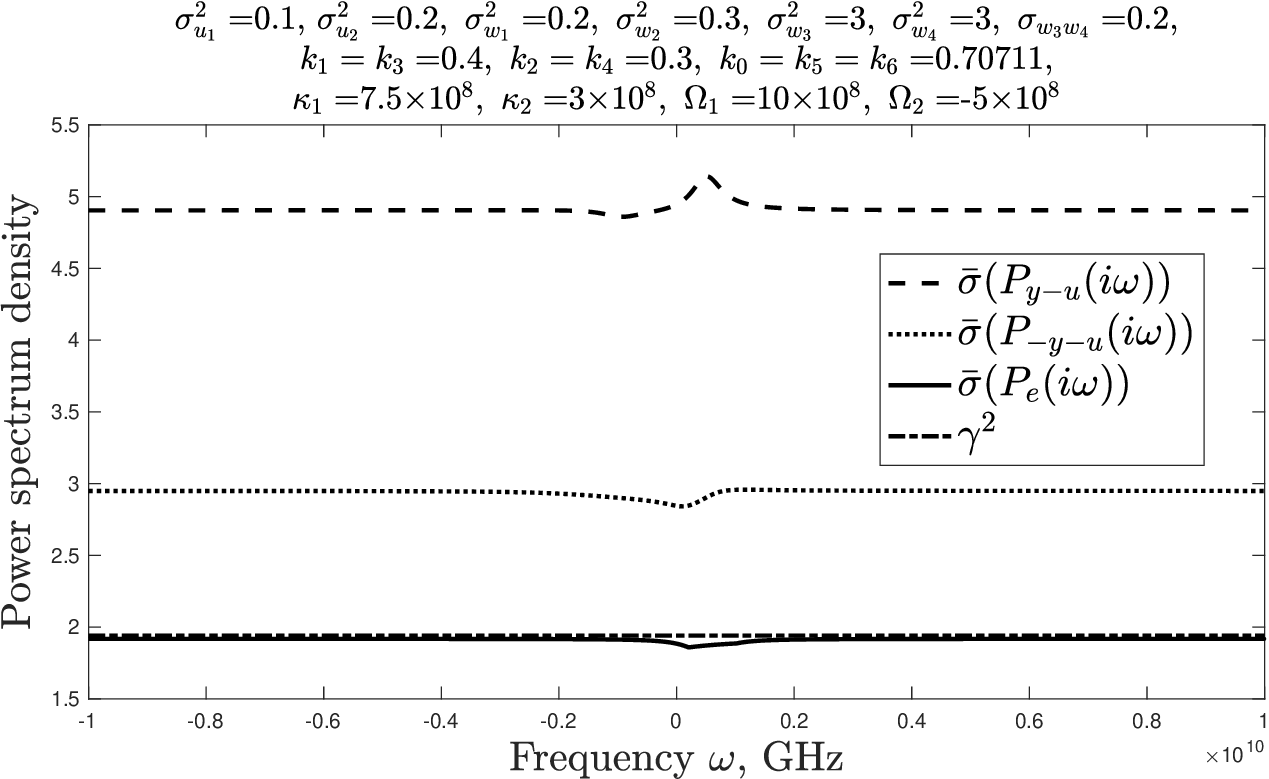}
  \caption{The largest eigenvalues of the PSD matrices $P_{y-u}$, $P_{-y-u}$ and
    $P_{e}$ and the suboptimal value $\gamma^2$ for the system in Fig.~\ref{two-cavity-system}.} 
  \label{fig:two-cavities-PSDs}
\end{figure}

\section{Conclusions}\label{Conclusions}

The paper has presented new results on the problem of coherent equalization for
completely passive quantum channels. It has explored benefits of  
restricting the class of candidate filters to those whose (1,1) block
$H_{11}(s)$ is contractive in the $H_\infty$ sense. We have shown that
doing so leads to near optimal physically realizable filters which
provide a tight upper bound on the mean-square optimal equalization
performance, to a desired accuracy.

Contractiveness has unlocked a host of tools
from control theory which we have been able to use in the design of
coherent equalizers. In particular, to obtain closed form formulas 
for components of the equalizing transfer function we have
employed results on spectral factorization of the transfer functions
$I-H_{11}(s)H_{11}(s)^H$ and $I-H_{11}(s)^HH_{11}(s)$~\cite{CALM-1997}. For
strictly contracting $H_{11}(s)$, the required spectral factors are easily
computed 
using stabilizing solutions to the algebraic Riccati equations~(\ref{eq:23})
and~(\ref{eq:24}). When 
$H_{11}(s)$ is not strictly contractive and is only 
nonexpanding in the sense of~(\ref{eq:1.UJ2}), such spectral factors are
substantially more difficult to compute~\cite{CALM-1997}. 
One of the possible directions for future research may be to investigate this
approach in a greater detail.

We have also shown that under the spectral factorization
Assumption~\ref{A.rho}, the design of a coherent quantum equalizer
reduces to finding a near optimal solution to a finite-dimensional
optimization problem. When the order of the filter is not restricted, the
algorithm boils down to solving a convex semidefinite program. This makes
the coherent equalizer design problem amenable to a wide range of
numerical methods 
for convex optimization. We have also been able to
dispense with some of the limitations of the previous approaches regarding
the dimensions of the filter input. 

Some of our results rely on the sufficient condition of
Theorem~\ref{SDP.primal.LMI}. We have
argued that this sufficient condition is appropriate in low
signal-to-noise ratio scenarios. The examples from quantum optics presented in
Section~\ref{examples} illustrate practical 
applicability of the proposed method in such situations. They
demonstrate that this condition is indeed satisfied when the system noise
has a high intensity, in which case a coherent filter could be found which
approximated the input field $u$ with a better mean-square accuracy than
the channel output 
$y$ does. Computing the exact threshold value on the signal-to-noise ratio below
which such improvement is achieved remains an open problem which is left for
future research.


\begin{thebibliography}{10}

\bibitem{Anderson-1967}
B.~D.~O. Anderson.
\newblock An algebraic solution to the spectral factorization problem.
\newblock {\em IEEE Transactions on Automatic Control}, 12(4):410--414, 1967.

\bibitem{BGR-2013}
J.~Ball, I.~Gohberg, and L.~Rodman.
\newblock {\em Interpolation of rational matrix functions}.
\newblock Birkh{\"a}user, Basel, 1990.

\bibitem{BCP-2021}
D.~Bertsimas, R.~Cory-Wright, and J.~Pauphilet.
\newblock Mixed-projection conic optimization: {A} new paradigm for modeling
  rank constraints.
\newblock {\em Operations Research}, 70(6):3321--3344, 2022.

\bibitem{LMI}
S.~Boyd, L.~El Ghaoui, E.~Feron, and V.~Balakrishnan.
\newblock {\em Linear Matrix Inequalities in System and Control Theory},
  volume~15 of {\em SIAM Studies in applied mathematics}.
\newblock SIAM, 1994.

\bibitem{CALM-1997}
D.~J. Clements, B.~D.~O. Anderson, A.~J. Laub, and J.~B. Matson.
\newblock Spectral factorization with imaginary-axis zeros.
\newblock {\em Linear Algebra and its Applications}, 250:225--252, 1997.

\bibitem{DP-2000}
G.~E. Dullerud and F.~Paganini.
\newblock {\em A Course in Robust Control Theory: {A} Convex Approach},
  volume~36 of {\em Texts in Applied Mathematics}.
\newblock Springer-Verlag, NY, 2000.

\bibitem{GJN-2010}
J.~E. Gough, M.~R. James, and H.~I. Nurdin.
\newblock Squeezing components in linear quantum feedback networks.
\newblock {\em Physical Review A}, 81(2):023804, 2010.

\bibitem{GZ-2015}
J.~E. Gough and G.~Zhang.
\newblock On realization theory of quantum linear systems.
\newblock {\em Automatica}, 59:139--151, 2015.

\bibitem{HSK-1999}
B.~Hassibi, A.~H. Sayed, and T.~Kailath.
\newblock {\em Indefinite-quadratic estimation and control: {A} unified
  approach to {$H^2$} and {$H^\infty$} theories}.
\newblock SIAM, Philadelphia, 1999.

\bibitem{HZ-2005}
R.~A. Horn and F.~Zhang.
\newblock Basic properties of the {Schur} complement.
\newblock In {\em The Schur Complement and Its Applications}, pages 17--46.
  Springer, 2005.

\bibitem{JNP-2008}
M.~R. James, H.~I. Nurdin, and I.~R. Petersen.
\newblock ${H}^\infty$ control of linear quantum stochastic systems.
\newblock {\em IEEE Transactions on Automatic Control}, 53(8):1787--1803, 2008.

\bibitem{Kailath-1981}
T.~Kailath.
\newblock {\em Lectures on Wiener and Kalman filtering}.
\newblock Springer, 1981.

\bibitem{KSH-2000}
T.~Kailath, A.~H. Sayed, and B.~Hassibi.
\newblock {\em Linear estimation}.
\newblock Prentice Hall, Upper Saddle River, NJ, 2000.

\bibitem{MP-2011}
A.~I. Maalouf and I.~R. Petersen.
\newblock Bounded real properties for a class of annihilation-operator linear
  quantum systems.
\newblock {\em IEEE Transactions on Automatic Control}, 56(4):786--801, 2011.

\bibitem{MJP-2016}
Z.~Miao, M.~R. James, and I.~R. Petersen.
\newblock Coherent observers for linear quantum stochastic systems.
\newblock {\em Automatica}, 71:264--271, 2016.

\bibitem{Nesterov-Nemirovskii}
Y.~E. Nesterov and A.~Nemirovsky.
\newblock {\em Interior Point Polynomial Methods in Convex Programming}.
\newblock SIAM, Philadelphia, PA, 1994.

\bibitem{Nurdin-2010}
H.~I. Nurdin.
\newblock On synthesis of linear quantum stochastic systems by pure cascading.
\newblock {\em IEEE Transactions on Automatic Control}, 55(10):2439--2444,
  2010.

\bibitem{NJD-2009}
H.~I. Nurdin, M.~R. James, and A.~C. Doherty.
\newblock Network synthesis of linear dynamical quantum stochastic systems.
\newblock {\em SIAM Journal on Control and Optimization}, 48(4):2686--2718,
  2009.

\bibitem{NY-2017}
H.~I. Nurdin and N.~Yamamoto.
\newblock {\em Linear Dynamical Quantum Systems}.
\newblock Springer, 2017.

\bibitem{Parthasarathy-2012}
K.~R. Parthasarathy.
\newblock {\em An introduction to quantum stochastic calculus}.
\newblock Birkh{\"a}user, 2012.

\bibitem{Rantzer-1996}
A.~Rantzer.
\newblock On the {K}alman--{Y}akubovich--{P}opov lemma.
\newblock {\em Systems \& Control Letters}, 28(1):7 -- 10, 1996.

\bibitem{SP-2012}
A.~J. Shaiju and I.~R. Petersen.
\newblock A frequency domain condition for the physical realizability of linear
  quantum systems.
\newblock {\em IEEE Transactions on Automatic Control}, 57(8):2033--2044, 2012.

\bibitem{Shaked-1990}
U.~Shaked.
\newblock {$H_\infty$}-minimum error state estimation of linear stationary
  processes.
\newblock {\em IEEE Transactions on Automatic Control}, 35(5):554--558, 1990.

\bibitem{UJ2b}
V.~Ugrinovskii and M.~R. James.
\newblock Active versus passive coherent equalization of passive linear quantum
  systems.
\newblock In {\em Proc. 58th IEEE CDC}, Nice, France, December 2019.
\newblock arXiv:1910.06462.

\bibitem{UJ2a}
V.~Ugrinovskii and M.~R. James.
\newblock Wiener filtering for passive linear quantum systems.
\newblock In {\em American Control Conference}, Philadelphia, PA, July 10-12
  2019.
\newblock arXiv:1901.09494.

\bibitem{UJ2}
V.~Ugrinovskii and M.~R. James.
\newblock Coherent equalization of linear quantum systems.
\newblock {\em Automatica}, 2022.
\newblock (submitted, under review). arXiv:2211.06003.

\bibitem{VP-2013}
I.~G. Vladimirov and I.~R. Petersen.
\newblock Coherent quantum filtering for physically realizable linear quantum
  plants.
\newblock In {\em 2013 European Control Conference (ECC)}, pages 2717--2723,
  2013.

\bibitem{VA-2014}
S.~L. Vuglar and H.~Amini.
\newblock Design of coherent quantum observers for linear quantum systems.
\newblock {\em New Journal of Physics}, 16(12):125005, 2014.

\bibitem{WM-2008}
D.~F. Walls and G.~J. Milburn.
\newblock {\em Quantum optics}.
\newblock Springer, 2008.

\bibitem{Wiener-1949}
N.~Wiener.
\newblock {\em The extrapolation, interpolation, and smoothing of stationary
  time series}.
\newblock Wiley, New York, 1949.

\bibitem{WM-2009}
H.~M. Wiseman and G.~J. Milburn.
\newblock {\em Quantum Measurement and Control}.
\newblock Cambridge University Press, 2009.

\bibitem{Yakubovich-1974}
V.~A. Yakubovich.
\newblock A frequency theorem for the case in which the state and control
  spaces are {Hilbert} spaces, with an application to some problems in the
  synthesis of optimal controls. {I}.
\newblock {\em Siberian Mathematical Journal}, 15(3):457--476, 1974.

\bibitem{Youla-1961}
D.~Youla.
\newblock On the factorization of rational matrices.
\newblock {\em IRE Transactions on Information Theory}, 7(3):172--189, 1961.

\bibitem{ZJ-2013}
G.~Zhang and M.~R James.
\newblock On the response of quantum linear systems to single photon input
  fields.
\newblock {\em IEEE Transactions on Automatic Control}, 58(5):1221--1235, 2013.

\end{thebibliography}

\newcommand{\noopsort}[1]{} \newcommand{\printfirst}[2]{#1}
  \newcommand{\singleletter}[1]{#1} \newcommand{\switchargs}[2]{#2#1}

\section*{Appendix}

This section presents the transfer function $H(s)$ computed using 
Algorithm~\ref{Alg1} for the quantum system in
Section~\ref{example2}.

The elements of the transfer function $H_{11}(s)$:
{\scriptsize\begin{eqnarray*}
&&  H_{11}^{(1,1)}(s)=  (-0.25247+5.5183e-14i)
\frac{(s+(4.499e08+1e09i))}
{(s+(6.847e08+1e09i))}, \\
&&  H_{11}^{(1,2)}(s)=  (0.25247+5.5247e-14i)
\frac{(s+(4.499e08+1e09i))}
{(s+(6.847e08+1e09i))}, \\ 
&&  H_{11}^{(2,1)}(s)=(-0.28664-3.3188e-14i)
\frac{(s+(2.647e08-5e08i))}
{(s+(3.147e08-5e08i))}, \\
&&  H_{11}^{(2,2)}(s)=  (-0.28664+3.3154e-14i)
\frac{(s+(2.647e08-5e08i))}
{(s+(3.147e08-5e08i))}. 
\end{eqnarray*}}

The elements of the transfer function $H_{12}(s)$:
{\scriptsize\begin{eqnarray*}
&&  H_{12}^{(1,1)}(s)=   -0.93049
\frac{(s+(7.126e08+1e09i))}
{(s+(6.847e08+1e09i))}, \\
&&  H_{12}^{(1,2)}(s)= (-1.0615e-14-2.6189e-14i)\\
&&\phantom{H_{21}}\times
\frac{(s+(2.085e09+4.247e08i))}
{(s+(6.847e08+1e09i))}, \\ 
&&  H_{12}^{(2,1)}(s)=(-1.0615e-14+2.6189e-14i)\\
&&\phantom{H_{21}}\times
\frac{(s+(2.288e09+1.922e08i))}
{(s+(3.147e08-5e08i))}, \\
&&  H_{12}^{(2,2)}(s)=-0.91415
\frac{(s+(3.236e08-5e08i))}
{(s+(3.147e08-5e08i))}. 
\end{eqnarray*}}

The elements of the transfer function $H_{21}(s)$:
{\scriptsize\begin{eqnarray*}
&&  H_{21}^{(1,1)}(s)=  0.92412 \\
&&\phantom{H_{21}}\times
\frac{(s+(3.191e08-5e08i)) (s-(3.236e08+5e08i))}{(s+(3.147e08-5e08i)) (s+(3.236e08-5e08i))}
\\
&&\phantom{H_{21}}\times
\frac{(s+(6.988e08+1e09i))(s-(7.126e08-1e09i))}
{(s+(6.847e08+1e09i))(s+(7.126e08+1e09i))}, \\
&&  H_{21}^{(1,2)}(s)= (-0.0099672+2.5374e-14i) \\
&&\phantom{H_{21}}\times
\frac{(s-(3.236e08+5e08i)) (s+(6.65e07-6.371e08i))}{(s+(3.147e08-5e08i)) (s+(3.236e08-5e08i))}
\\
&&\phantom{H_{21}}\times
\frac{(s-(7.126e08-1e09i)) (s+(1.829e09+1.137e09i))}
{(s+(6.847e08+1e09i))(s+(7.126e08+1e09i))},\\
&&  H_{21}^{(2,1)}(s)=(-0.0099672-2.5374e-14i)\\
&&\phantom{H_{21}}\times
\frac{(s-(3.236e08+5e08i)) (s+(6.65e07-6.371e08i))}{(s+(3.147e08-5e08i)) (s+(3.236e08-5e08i))}
\\
&&\phantom{H_{21}}\times
\frac{(s-(7.126e08-1e09i)) (s+(1.829e09+1.137e09i))}
{(s+(6.847e08+1e09i))(s+(7.126e08+1e09i)) }, \\
&&  H_{21}^{(2,2)}(s)=0.92412 \\
&&\phantom{H_{21}}\times
\frac{(s+(3.191e08-5e08i)) (s-(3.236e08+5e08i))}{(s+(3.147e08-5e08i)) (s+(3.236e08-5e08i))}
\\
&&\phantom{H_{21}}\times
\frac{(s+(6.988e08+1e09i))(s-(7.126e08-1e09i))}
{(s+(6.847e08+1e09i))(s+(7.126e08+1e09i)) }. 
\end{eqnarray*}}

The elements of the transfer function $H_{22}(s)$:
{\scriptsize\begin{eqnarray*}
&&  H_{22}^{(1,1)}(s)= (-0.25247-5.5183e-14i) \\
&&\phantom{H_{22}}\times
\frac{(s-(3.236e08+5e08i)) (s-(4.499e08-1e09i))}
{(s+(3.236e08-5e08i)) (s+(6.847e08+1e09i))}, \\
&&  H_{22}^{(1,2)}(s)= (-0.28664+3.3188e-14i)\\
&&\phantom{H_{22}}\times 
\frac{(s-(2.647e08+5e08i)) (s-(7.126e08-1e09i))}{(s+(3.147e08-5e08i)) (s+(7.126e08+1e09i))},
\\
&&  H_{22}^{(2,1)}(s)=(0.25247-5.5247e-14i)\\
&&\phantom{H_{22}}\times
\frac{(s-(3.236e08+5e08i)) (s-(4.499e08-1e09i))}{(s+(3.236e08-5e08i)) (s+(6.847e08+1e09i))},\\
&&  H_{22}^{(2,2)}(s)=(-0.28664-3.3154e-14i) \\
&&\phantom{H_{22}}\times
\frac{(s-(2.647e08+5e08i)) (s-(7.126e08-1e09i))}{(s+(3.147e08-5e08i)) (s+(7.126e08+1e09i))}. 
\end{eqnarray*}}

\end{document}